\documentclass[10pt]{amsart}

\usepackage{xr}
\usepackage{amssymb}
\usepackage{mathdefs}
\usepackage{mathtools}
\usepackage{algorithm}
\usepackage{algpseudocode}
\usepackage{subcaption}
\usepackage{xcolor}

\newcommand{\Caratheodory}{Carath\'{e}odory}
\newcommand{\nnls}{\texttt{NNLS}}
\newcommand{\lp}{\texttt{LP}}
\newcommand{\csp}{\texttt{CSP}}
\newcommand{\scsp}{\texttt{SCSP}}
\newcommand{\gscsp}{\texttt{GSCSP}}

\usepackage[backend=biber,date=year,giveninits=true,sorting=nyt,natbib=true,maxcitenames=2,maxbibnames=10,url=false,doi=true,backref=false,sortcites=true]{biblatex}
\title[Efficient and robust Carath\'{e}odory-Steinitz pruning]{Efficient and robust Carath\'{e}odory-Steinitz pruning of positive discrete measures}
\author{Filip B\v{e}l\'ik
\and {Jesse Chan}
\and {Akil Narayan}
}
\usepackage{hyperref}
\usepackage{cleveref}
\usepackage{booktabs}

\newtheorem{theorem}{Theorem}[section]

\newtheorem{definition}{Definition}[section]

\theoremstyle{remark}
\newtheorem{remark}{Remark}[section]
\newtheorem{assumption}{Assumption}[section]

\newtheoremstyle{example}{1em plus .2em minus .1em}{1em plus .2em minus .1em}%
 {}%         Body font
 {}%         Indent amount (empty = no indent, \parindent = para indent)
 {\bfseries}% Thm head font
 {}%        Punctuation after thm head
 {\newline}%     Space after thm head (\newline = linebreak)
 {\thmname{#1}\thmnumber{ #2}\thmnote{ #3}}%         Thm head spec
\begin{document}
\maketitle
\begin{abstract}
    In many applications, one seeks to approximate integration against a positive measure of interest by a positive discrete measure: a numerical quadrature rule with positive weights. One common desired discretization property is moment preservation over a finite dimensional function space, e.g., bounded-degree polynomials. Carath\'{e}odory's theorem asserts that if there is any finitely supported quadrature rule with more nodes than the dimension of the given function space, one can form a smaller (and hence more efficient) positive, nested, quadrature rule that preserves the moments of the original rule. 
    
We describe an efficient streaming procedure for Carath\'{e}odory-Steinitz pruning, a numerical procedure that implements Carath\'{e}odory's theorem for this measure compression. The new algorithm makes use of Givens rotations and on-demand storage of arrays to successfully prune very large rules whose storage complexity only depends on the dimension of the function space. This approach improves on a naive implementation of Carath\'{e}odory-Steinitz pruning whose runtime and storage complexity are quadratic and linear, respectively, in the size of the original measure. We additionally prove mathematical stability properties of our method with respect to a set of admissible, total-variation perturbations of the original measure. Our method is compared to two alternative approaches with larger storage requirements: non-negative least squares and linear programming, and we demonstrate comparable runtimes, with improved stability and storage robustness. Finally, we demonstrate practical usage of this algorithm to generate quadrature for discontinuous Galerkin finite element simulations on cut-cell meshes.

\end{abstract}

\section{Introduction}
\label{sec:intro}
Efficient and accurate quadrature (or cubature) rules that approximate integrals are fundamental ingredients in computational science, being used for numerical or statistical integration in the context of solutions of differential equations, uncertainty quantification, inference, and scientific machine learning. In these application scenarios one may have access to an acceptably-accurate quadrature rule with positive weights; the challenge is that this quadrature rule might be too large to use in practice because it requires too many function evaluations. To ameliorate this situation, one can consider using this starting quadrature rule to identify a quadrature rule with many fewer nodes that retains desirable properties, in particular retains both positivity and accuracy, where the latter is quantified by exact integration of specified moments. The core algorithm we consider, \Caratheodory-Steinitz pruning (CSP), is one strategy that identifies a quadrature rule that is nested with respect to the original (hence, is a ``pruned'' version because nodes are removed) \cite{steinitz_bedingt_1913}. The CSP algorithm has been particularly popular for its clear and easy implementation, and has seen applications in contexts requiring high-dimensional quadrature over general domains \cite{davis_construction_1967,wilson_general_1969,piazzon_caratheodorytchakaloff_2017,tchernychova_caratheodory_2016,vandenbos_nonintrusive_2017,vandenbos_generating_2020,glaubitz_stable_2020-1,glaubitz_constructing_2020,glaubitz_construction_2023,GONOSKOV2022108200,Sommariva02092015,Elefante2022CQMCAI,feldman_2019,statman_2023}.

However, a primary challenge with the CSP algorithm is computational cost. If an $M$-point positive quadrature rule is pruned to an $N$-point positive quadrature rule subject to $N$ moment constraints, then a naive implementation requires a cumulative $\mathcal{O}((M-N) M N^2)$ computational complexity and $\mathcal{O}(M N)$ storage complexity. In several practical use cases of interest, $M \gg N$, which makes both the storage and complexity demands of a naive CSP algorithm onerous.

Our contributions in this paper are the following two major advances: First, we devise a compute- and storage-efficient version of CSP, which makes the per-step computational complexity independent of $M$ and improves overall storage requirements to $\mathcal{O}(N^2)$ when the algorithm is used in streaming contexts for pruning a size-$M$ positive quadrature rule down to $N$ nodes.
Our storage-efficient, ``streaming'' version of CSP, the \scsp{} algorithm, is given in \Cref{alg:scsp}. A similar streaming idea has been used previously for the generation of so-called implicit quadrature rules \cite{vandenbos_generating_2020} and in forming least-mean-squares coresets \cite{statman_2023}. A further augmentation of this algorithm, the \gscsp{} procedure (``Givens \scsp{}''), is an efficient procedure for computing cokernel vectors. The \gscsp{} algorithm requires only $\mathcal{O}(N^2)$ complexity per iteration for a cumulative $\mathcal{O}((M-N) N^2) + \mathcal{O}(N^3)$ computational complexity. This efficiency is gained by exercising Givens rotations for updating cokernel vectors of a matrix. The \gscsp{} algorithm is \Cref{alg:scsp} with the augmentation in Supplementary Algorithm S1.

Second, we provide a new stability guarantee for the \scsp{} and \gscsp{} algorithms: by considering any particular quadrature rule as a (discrete) measure, we show that these procedures are mathematically stable with respect to small total variation perturbations of the input measure, namely perturbations of the existing weights and the insertion of new nodes with small weights, provided the ordering of the original nodes is preserved. Under these perturbations, both algorithms are locally Lipschitz (and in particular continuous) with respect to the total variation distance on both the input and the output. See \Cref{thm:muM-perturbation}.

In the numerical results presented in \Cref{sec:numerical}, we demonstrate that the \gscsp{} algorithm can successfully prune very large rules with one billion points and compare the computational efficiency of \gscsp{} to competing algorithms, in particular, a non-negative least squares (NNLS) formulation and a linear programming (LP) formulation. We note that modifications to the ``standard'' NNLS approach have also been proposed that decrease runtime by up to one order of magnitude \cite{Elefante2022CQMCAI}. We also provide supporting evidence for the total variation stability guarantee for \scsp{} and \gscsp, and show that the stability properties of this algorithm are more favorable than the stability properties of the alternative NNLS and LP algorithms. We demonstrate the potential of our new algorithm by generating nontrivial quadrature on two-dimensional cut-cell finite element mesh geometries. The \gscsp{} and other related ``pruning'' algorithms are implemented in the open-source software package \texttt{CaratheodoryPruning.jl}.

\section{Background}
We use the notation $\N \coloneqq \{1, 2, \ldots, \}$ and $\N_0 \coloneqq \{0\} \cup \N$. For $N \in \N$, we let $[N] = \{1, \ldots, N\}$. Lowercase/uppercase boldface letters are vectors/matrices, respectively, e.g., $\bs{x}$ is a vector and $\bs{X}$ is a matrix. If $\bs{A} \in \R^{M \times N}$ with $S \subset [M]$ and $T \subset [N]$, then we use the notation,
\begin{align*}
  \bs{A}_{S*} &\in \R^{|S| \times N}, & 
  \bs{A}_{*T} &\in \R^{M \times |T|}, &
  \bs{A}_{ST} &\in \R^{|S| \times |T|},
\end{align*}
to slice $\bs{A}$ producing, respectively, the $S$-indexed rows of $\bs{A}$, the $T$-indexed columns of $\bs{A}$, and the submatrix formed by rows indexed $S$ and columns indexed $T$. Throughout, we will consider $S$, $T$, and any other subsets of indices as ordered sets (e.g., sequences) so that, e.g., the first row of $\bs{A}_{S\ast}$ is the row of $\bs{A}$ corresponding to the first index in the ordered set $S$. Similarly, given a vector $\bs{b}\in\R^M$, we use the notation $\bs{b}_S\in\R^{|S|}$ to slice $\bs{b}$, producing the ordered, $S$-indexed, elements of $\bs{b}$.  If $\bs{v}$ and $\bs{w}$ are vectors of the same size, then $\bs{v} \geq \bs{w}$ means that the inequality holds component-wise. Unless otherwise noted, we will denote the two-norm of a vector as $\|\bs{w}\| = \|\bs{w}\|_2 = \sqrt{\bs{w}^T\bs{w}}$. We denote the nonnegative reals by $\R_+$ and the positive reals by $\R_{++}$.  

Given a set of points $P = \{\bs{p}_1,\bs{p}_2,\ldots,\bs{p}_M\} \subset \R^N$, we say that a point $\bs{p}\in\R^N$ lies in the conic hull of $P$, denoted $\bs{p}\in\text{cone}(P)$, if there exist $\{w_m\}_{m\in[M]}\subset\R_+$ such that
\[\bs{p} = \sum_{m\in[M]} w_m \bs{p}_m.\]

\subsection{Positive quadrature rules: Tchakaloff's theorem}\label{ssec:tchakaloff}
Let $(X,\mathcal{M},\mu)$ be a measure space, $\mu$ is a positive measure, e.g., $\mu$ a probability measure, and let $V$ be an $N$-dimensional subspace of functions in $L_\mu^1(X)$ spanned by basis elements $v_j$:
\begin{align}\label{eq:V-def}
  V &\coloneqq \mathrm{span} \{v_1, \ldots, v_N \}, & v_j &: X \rightarrow \R.
\end{align}
Our main goal is to construct a \textit{positive} quadrature rule, i.e., a set of nodes and weights, $\{x_q\}_{q \in [Q]} \subset X$ and $\{w_q\}_{q \in [Q]} \subset \R_{++}$, such that, 
\begin{align}\label{eq:tchakaloff}
  \int_X v(x) \dx{\mu}(x) = \sum_{q \in [Q]} w_q v(x_q), \hskip 15pt \forall\;\; v \in V.
\end{align}
The above procedure is sometimes called \textit{measure compression}.
We assume $V$ to be a subspace of $L_\mu^1(X)$ containing fixed representatives for which point evaluation is well-defined; a specialization would be a subspace of $\mu$-integrable, continuous functions.
The classical Tchakaloff's Theorem states that \eqref{eq:tchakaloff} is attainable for $Q \leq N$ for general polynomial subspaces $V$ \cite{tchakaloff_formules_1957,davis_construction_1967,putinar_note_1997,curto_duality_2002-1,bayer_proof_2006}. This result has also recently been extended to general subspaces $V$ \cite{schafer_tchakaloff_2025}. These results are generally not computationally constructive when $\mu$ is not a finitely supported measure.

The central problem statement of this paper is that we seek to computationally realize Tchakaloff measure compression for general subspaces $V$ but where $\mu$ is a finitely-supported discrete measure. In general, the bound $Q = N = \dim(V)$ is sharp, and so we focus on computationally attaining this goal in this paper. (However, $Q < N$ is attainable in some special cases, see Section 2 of the Supplementary Information for an example.)

We therefore first assume that some quadrature rule with more than $\dim(V)$ nodes is available that meets the accuracy requirements \eqref{eq:tchakaloff}. Equivalently, we make the fairly strong assumption that the initial measure $\mu$ is a finitely-supported (discrete) measure (or can be approximated sufficiently well by a finitely-supported measure), and seek to compress this measure subject to a $V$-moment matching condition. 

\subsection{Finitely supported measures}
In many scenarios one is able to construct a positive quadrature rule with $M \gg N$ points; another way to state this is that there is a measure $\mu_M$ supported on $M$ points, i.e., 
\begin{align}\label{eq:mu-finitely-supported}
  \mu &\approx \mu_M, & \dx{\mu_M}(x) &= \sum_{m \in [M]} w_m \delta_{x_m}, & M < \infty,
\end{align}
where $\delta_{x}$ is the Dirac mass centered at $x$ and $\{x_m\}_{m\in[M]} \subset X$ and $\{w_m\}_{m\in[M]} \subset \R_{++}$ are the nodes and weights for $\mu_M$ respectively. If $M \leq N$, we already have a Tchakaloff-realizing quadrature rule, so we assume without loss of generality that $M > N$. In this case we have that $\mu_M$, defined by its nodes and weights, has certain moments of $V$.
With a fixed $N$-dimensional subspace $V$ with basis $\{v_n\}_{n\in[N]}$, define,
\begin{align}\label{eq:eta_n}
  \eta_n \coloneqq \int_{X} v_n(x) \dx{\mu_M}(x) = \sum_{m \in [M]} w_m v_n(x_m) \in \R, \hskip 15pt n \in [N].
\end{align}
A Tchakaloff theorem for finitely supported measures states the following.
\begin{theorem}\label{thm:finite-tchakaloff}
  Let $(\mu_M,V)$ be as described above with finite moments as defined in \eqref{eq:eta_n}. Then \eqref{eq:tchakaloff} holds with $Q \leq N = \dim(V)$ where the $Q$ quadrature nodes are a subset of $\mathrm{supp} (\mu_M) = \{x_m\}_{m \in [M]}$.
\end{theorem}
The above result is not new and is essentially well-known. See, e.g, related statements in \cite{tchernychova_caratheodory_2016,piazzon_caratheodorytchakaloff_2017,steinitz_bedingt_1913}. 
One standard proof of \Cref{thm:finite-tchakaloff} reveals a popular algorithm that makes the result constructive; this proof relies on a minor variant of Carath\'{e}odory's theorem in convex geometry.
\begin{theorem}[Carath\'{e}odory's theorem, conic version \cite{eisenbrand_2006}]\label{thm:cara-conic}
  Let $P \subset \R^N$ be a finite set of points in $\R^N$ with $|P| > N$. If $\bs{p} \in \cone(P)$, then there exists a subset $S \subset P$ with $|S| \leq N$ such that $\bs{p} \in \cone(S)$.
\end{theorem}
\begin{remark}
The more traditional phrasing of Carath\'{e}odory's Theorem that considers the stronger notion of \textit{convex} combinations yields the looser conclusion $|S| \leq N+1$.
\end{remark}
To see how this applies to our situation, we provide a direct, simple proof that reveals a computational implementation. Like \Cref{thm:finite-tchakaloff} itself, neither this proof nor the algorithm are new. 

\subsection{The Carath\'{e}odory-Steinitz ``pruning'' construction}\label{ssec:cs-pruning}
In this section we review one simple constructive proof of both \Cref{thm:finite-tchakaloff} and \Cref{thm:cara-conic} revealing an algorithm. This algorithm has recently seen considerable use \cite{piazzon_caratheodorytchakaloff_2017,tchernychova_caratheodory_2016,glaubitz_stable_2020-1,glaubitz_constructing_2020}. We attribute this algorithm originally to Steinitz \cite{steinitz_bedingt_1913}, and will hence refer to the following naive algorithm as the \textit{Carath\'{e}odory-Steinitz pruning} (CSP) algorithm.

If $M \leq N$, then \Cref{thm:finite-tchakaloff} is trivially proven, so without loss we assume $M > N$. The core idea is the simple observation that the moment conditions \eqref{eq:eta_n} can be written through linear algebra:
\begin{align}\label{eq:moment-conditions}
  \bs{V}^T \bs{w} &= \bs{\eta}, & 
  \bs{V} &= \left(\begin{array}{ccc} 
    \horzbar & \bs{v}(x_1)^T & \horzbar \\
    \horzbar & \bs{v}(x_2)^T & \horzbar \\
             & \vdots &               \\
    \horzbar & \bs{v}(x_M)^T & \horzbar \\
  \end{array}\right) \in \R^{M \times N},
\end{align}
where $\bs{w}, \bs{v}(x_m)$, and $\bs{\eta}$ are
\begin{align}\label{eq:v-def}
  \bs{v}(x_m) &\coloneqq \left( v_1(x_m), \; \ldots, \; v_N(x_m) \right)^T \in \R^N, \hskip 5pt m \in [M], \\\nonumber
  \bs{w} &\coloneqq \left( w_1, \; \ldots, w_M \right)^T \in \R^M_{++}, \\\nonumber
  \bs{\eta} &\coloneqq \left( \eta_1, \; \ldots, \; \eta_N \right)^T \in \R^N,
\end{align}
with $\eta_n$, $n\in[N]$, defined in \eqref{eq:eta_n}. 
If $M > N$, then $\bs{V}^T \in \R^{N \times M}$ has a non-trivial kernel, so there is some kernel vector, say $\bs{n} \neq \bs{0}$, such that,
\begin{align*}
  \bs{V}^T \left( \bs{w} - c \bs{n} \right) &= \bs{\eta}, & \forall \; c &\in \R.
\end{align*}
This kernel vector can be used to construct a size-$(M-1)$ quadrature rule by augmenting $\bs{w}$. We first partition $[M]$ into sets where $\bs{n}$ is positive, negative, and 0,
  \begin{align*}
    S_{\pm} &\coloneqq \left\{ m \in [M] \;\big|\; \pm n_m > 0 \right\},  &
    S_0 &\coloneqq \left\{ m \in [M] \;\big|\; n_m = 0 \right\}, & 
    [M] &= S_+ \cup S_- \cup S_0.
  \end{align*}
Because $\bs{n} \neq \bs{0}$, it is not possible for both $S_+$ and $S_-$ to be empty. We now define smallest-magnitude constants $c$ that ensure $\bs{w} - c \bs{n}$ has (at least) one zero component:
\begin{align}\label{eq:mc-def}
  m_{\pm} &= \argmin_{m \in S_{\pm}} \left|\frac{w_m}{n_m}\right|, & 
  c_{\pm} &= \frac{w_{m_\pm}}{n_{m_\pm}},
\end{align}
where when $S_{\pm} = \emptyset$ we assign $c_{\pm} = \pm\infty$. With this construction, $c_- < 0 < c_+$, and,
\begin{align*}
  c \in (c_-, c_+) \hskip 5pt &\Longleftrightarrow \hskip 5pt \bs{V}^T (\bs{w} - c \bs{n}) = \bs{\eta} \textrm{ and } \bs{w} - c \bs{n} > \bs{0}, \\
  c = c_{\pm} \hskip 5pt &\Longleftrightarrow \hskip 5pt \bs{V}^T (\bs{w} - c \bs{n}) = \bs{\eta}, \; \bs{w} - c \bs{n} \geq \bs{0}, \; \textrm{ and } (\bs{w} - c \bs{n})_{m_\pm} = 0,
\end{align*}
where in the second line we assume that both $c_{\pm}$ are finite; at least one of them must be finite since $S_+ \cup S_-$ is non-empty. Hence, choosing either $c = c_+$ or $c = c_-$, setting $\bs{w} \gets \bs{w} - c \bs{n}$, and then removing row $m_{\pm}$ (and all other zeroed rows) from both $\bs{w}$ and $\bs{V}$, constructs an (at most) $(M-1)$-point rule with $\bs{w} \geq \bs{0}$ satisfying $\bs{V}^T\bs{w} = \bs{\eta}$. A visualization of the two pruning choices is provided in the supplementary materials. This process can be repeated while $\bs{V}$ has a nontrivial cokernel, which is generically until $\bs{V}$ is square, corresponding to an (at most) $Q=N$-point rule, completing the proofs of both \Cref{thm:finite-tchakaloff} and \Cref{thm:cara-conic}.

The sign $\sigma \in \{+, -\}$ that identifies $c_\sigma$ must be algorithmically chosen. In general, this choice is given by
\begin{align}\label{eq:sigma}
  \sigma &= \left\{\begin{array}{rl}
    +, \hspace{4mm}& S_- = \emptyset \\
    -, \hspace{4mm}& S_+ = \emptyset \\
    \texttt{SigSelect}(\bs{V}, \bs{w}, \bs{n}), \hspace{4mm}& \textrm{otherwise.}
  \end{array}\right.
\end{align}
One \textit{example} of the function \texttt{SigSelect} would be the simple rule,
\begin{align}\label{eq:sigselect-example}
  \texttt{SigSelect} = \argmin_{\sigma \in \{+,-\}} \left|c_{\sigma}\right|  \hskip 10pt \Longrightarrow \hskip 10pt m = \argmin_{m \in S_+ \cup S_-} \frac{w_m}{|n_m|}, \hskip 5pt c = \frac{w_m}{n_m},
\end{align}
which simply chooses $+$ or $-$ based on which choice corresponds to a minimum-norm perturbation of $\bs{w}$. In general, this choice could depend on $\bs{V}$, $\bs{w}$, and $\bs{n}$. Pseudocode for the CSP procedure is given in \Cref{alg:csp}.

\begin{algorithm}
  \caption{\texttt{CSP}: \Caratheodory-Steinitz pruning}\label{alg:csp}
    {\raggedright\noindent
        \textbf{Input: } $\bs{V} \in \mathbb{R}^{M\times N}$, $\bs{w} \in \mathbb{R}^M_{++}$\\
        \textbf{Output: } $S$ with $|S| \leq N$, $\bs{w}_S \in \R^{|S|}_{++}$\par}
    \begin{algorithmic}[1]
        \State $S = [M]$
        \While{$|S| > N$}
          \State Compute ${\bs{n}} \in \mathrm{ker}({\bs{V}_{S*}}^T)$. \Comment{$\mathcal{O}(|S| N^2)$}\label{lst:kernel-vector}
          \State Identify $S_\pm$ and compute $c_{\pm}$ in \eqref{eq:mc-def} using $S_{\pm}$, $\bs{w}_S$, $\bs{n}$.
          \State Choose $\sigma \in \{+,-\}$ as in \eqref{eq:sigma} \Comment{\texttt{SigSelect}, e.g., as in \eqref{eq:sigselect-example}}
          \State Set $\bs{w}_S \gets \bs{w}_S - c_\sigma \bs{n}$, $P = \left\{ s \in S \;\;\big|\;\; w_s = 0 \right\}$
          \State $S \gets S \backslash P$
        \EndWhile
    \end{algorithmic}
\end{algorithm}

\begin{remark}
  One can continue the while loop in \Cref{alg:csp} with $|S| \leq N$ so long as $\bs{V}^T$ has a non-trivial kernel. This would yield a rule with $|S| < N$ points. However, if a positive quadrature rule of size $|S| < N$ does exist, there is no guarantee that \Cref{alg:csp} finds this rule, and instead it can terminate with $N$ nodes.
\end{remark}

A direct implementation of \Cref{alg:csp} requires $\mathcal{O}(M N)$ storage, largely to access the full original matrix $\bs{V}$. The computational complexity is $\mathcal{O}\left(M (M-N) N^2\right) \allowbreak \lesssim \mathcal{O}(M N^3 + M^2 N^2)$, since the dominant cost is identification of the kernel vector $\bs{n}$ at each step. Note that the most expensive step is when $|S|=M$ and that the algorithm terminates in a maximum of $(M-N)$ steps. This algorithm has been extensively explored in \cite{vandenbos_generating_2020}, which provides accompanying software.

In contrast, the main algorithmic innovation of this paper is a procedure that accomplishes the same result as the \csp{} algorithm but requires only $\mathcal{O}(N^2)$ storage and $\mathcal{O}((M-N) N^2) + \mathcal{O}(N^3) \lesssim \mathcal{O}(MN^2 + N^3)$ complexity. In particular, the new algorithm has a storage complexity independent of $M$ and a computational complexity that is linear in $M$, which is of considerable benefit in the realistic $M \gg N$ setting.

\subsection{Alternative algorithms}\label{ssec:alternatives}
We describe two alternatives to CSP that have also enjoyed popularity due to their computational convenience \cite{piazzon_caratheodorytchakaloff_2017,tchernychova_caratheodory_2016,glaubitz_stable_2020-1,glaubitz_constructing_2020}.

The first alternative method employs non-negative least squares (NNLS), and numerically solves the quadratic programming problem,
\begin{align}\label{eq:nnls}
  \argmin_{\bs{v}\in\R^{M}_+} \left\| \bs{V}^T \bs{v} - \bs{\eta} \right\|^2.
\end{align}
In order to accomplish pruning, one hopes that $\bs{v}$ is $N$-sparse. The explicit optimization formulation above does not suggest why such sparse solutions should be obtained, but practical algorithms to solve this problem implicitly enforce sparsity \cite{lawson_1995}, and in generic situations numerically solving \eqref{eq:nnls} indeed produces $N$-sparse solutions and achieves zero objective. There are some optimizations of NNLS-type algorithms that can accelerate computations \cite{bro1997fast,dessole_lawson_hanson}. Another NNLS approach that resembles the streaming improvement that will be introduced in the next section makes use of a sequence of nested, increasing cardinality, subsets of the quadrature nodes \cite{Elefante2022CQMCAI}. The approach is promising and is found to obtain up to one order of magnitude runtime speedup compared to standard NNLS in most examples, requiring only a few iterations \cite{Elefante2022CQMCAI}. However, the algorithm is dependent on a node enrichment procedure which may affect the number of iterations required and the resulting runtime.

A second alternative is through linear programming. 
Linear programming has previously found application in computing Gauss-type quadrature in \cite{ryu_extensions_2015}. 
First, observe that any solution to the linear moment constrained problem \eqref{eq:moment-conditions} with the desired non-negativity constraints is given by,
\begin{subequations}
\begin{align}
  \label{eq:W1}
  W \coloneqq& \left\{ \bs{v} \in \R^M_+ \;\;\big|\;\; \bs{V}^T \bs{v} = \bs{\eta}\right\} \\\label{eq:W2}
    =& \left\{ \bs{w} + \bs{K} \bs{z} \in \R^M_+ \;\;\big|\;\; \bs{z} \textrm{ arbitrary}\right\},
\end{align}
\end{subequations}
where $\bs{K}$ is a(ny) matrix whose range is the cokernel of $\bs{V}$, and $\bs{w}$ is a(ny) set of $M$ weights that matches moments. Hence, the feasible set $W$ of weight vectors $\bs{v}$ is a polyhedron in $\R^M$ of dimension $\dim (\coker (\bs{V})) \geq M-N$. The extreme points of this polyhedron correspond to at least $M-N$ active constraints in $\R^M_+$, i.e., at least $M-N$ zero weights. Hence, one way to identify a vector of quadrature weights with at most $N$ entries is to identify one point in $\mathrm{ex}(W)$, the set of extreme points of $W$, which can be accomplished through linear programming: for some $\bs{c} \in \R^M$, solving,
\begin{align}\label{eq:lp}
  \min_{\bs{v}} \bs{c}^T \bs{v} \textrm{ subject to } \bs{v} \in W,
\end{align}
generically produces an $N$-sparse solution. ``Generically'' means, e.g., that if $\bs{c}$ has random components that are drawn iid from a standard normal distribution, then with probability 1 the solution to \eqref{eq:lp} has at most $N$ non-zero entries. Note that \eqref{eq:lp} does not require $\bs{K}$ through definition \eqref{eq:W1} of $W$, but having knowledge of $\bs{K}$ converts \eqref{eq:lp} from an $M$-dimensional optimization to a $(M-N)$-dimensional one on the variables $\bs{z}$ through definition \eqref{eq:W2} of $W$. A similar linear programming formulation is used in the integration of parameter-dependent functions in \cite{yano_lp}.

One method that builds on top of an existing pruning algorithm is called Fast-Carath\'{e}odory-set \cite{feldman_2019}. Fixing an integer $k\approx e N$ (where $e$ is Euler's number), the algorithm works by partitioning the indices $[M]$ into $k$ disjoint subsets of near equal size, $\{P_i\}_{i=1}^k$. Upon forming weighted means of each cluster, pruning can be performed on the size $(k,N)$ problem
\[
  \widetilde{\bs V}^T\bs 1 = 
  \sum_{i=1}^k 1\cdot \widetilde{\bs v}_i = \bs\eta,\quad \widetilde{\bs v}_i = \sum_{j \in P_i} w_j \bs v_j.
\]
After pruning down to $N$ clusters and expanding the original terms, the size of the problem has reduced from $M$ to about $MN/k$. Repeating, the problem sizes reduce by a factor of $k/N$ after each iteration. Forming the weighted means has a total runtime complexity of $\mathcal{O}(MN)$. Letting $t(M,N)$ denote the runtime complexity of pruning a size $(M,N)$ problem, the pruning of size $(k,N)$ problems is repeated about $\log_{k/N}(M)=\log(M)/\log(k/N)$ times. Altogether, taking $k\approx eN$ results in a runtime complexity of $\mathcal{O}(MN + t(eN,N)\log(M))$. A naive implementation yields $t(eN,N)=\mathcal{O}(N^4)$; incorporating the proposed Givens procedure described next would reduce this complexity to $\mathcal{O}(N^3)$.

\section{Kernel vector computations through Givens rotations}\label{sec:csp-algorithms}
We present the main algorithmic novelty of this paper in this section, which is an efficient procedure to accomplish line \ref{lst:kernel-vector} in \Cref{alg:csp}, i.e., to compute cokernel vectors repeatedly for progressively row-pruned matrices $\bs{V}$. One essential idea is that computing cokernel vectors is equivalent to computing vectors orthogonal to the range, and the latter is accomplished through the QR decomposition of a row rank-deficient matrix. For an $M \times N$ matrix $\bs{V}$ with $M > N$ and rank $N$:
\begin{align*}
  \bs{V} &= \bs{Q} \bs{R} = \left[ \bs{Q}_1 \; \bs{Q}_2 \right] \bs{R}, & \bs{Q}_1 \in \R^{M \times N}, \;\; \bs{Q}_2 \in \R^{M \times (M-N)}, \;\; \bs{R} \in \R^{M \times N},
\end{align*}
where $\bs{Q}$ has orthonormal columns and $\bs{R}$ is upper triangular. The matrix $\bs{Q}_2$ above and the matrix $\bs{K}$ in \eqref{eq:W2} have the same range; the difference is that $\bs{Q}_2$ has orthonormal columns. A(ny) nontrivial vector in the range of $\bs{Q}_2$ is a kernel vector of $\bs{V}^T$, equivalently is a cokernel vector of $\bs{V}$. Hence, we can identify a cokernel vector through a QR decomposition of $\bs{V}$, having complexity $\mathcal{O}(M N^2)$.

\subsection{The \scsp{} algorithm: \texorpdfstring{$\mathcal{O}(N^2)$}{O(N\^{2})} storage}
One simple modification of the above approach is motivated by observing that it is wasteful to compute $M-N$ kernel vectors (all of $\bs{Q}_2$) when only 1 is needed. Even if only a single cokernel vector is computed rather than all of $\bs{Q}_2$, the QR decomposition of the full matrix $\bs{V}$ still requires storage complexity that depends on $M$. An alternative and more efficient approach is to compute a single cokernel vector for a slicing of $\bs{V}$. If $S \subset [M]$ with $|S| > N$ is any row subset, then consider the QR decomposition of the $S$-row sketched matrix $\bs{V}_{S\ast}$:
\begin{align*}
  \bs{V}_{S\ast} &= \widetilde{\bs{Q}} \widetilde{\bs{R}} = \left[ \widetilde{\bs{Q}}_1 \; \widetilde{\bs{Q}}_2 \right] \widetilde{\bs{R}}, & \widetilde{\bs{Q}}_1 \in \R^{|S| \times N}, \;\; \widetilde{\bs{Q}}_2 &\in \R^{|S| \times (|S|-N)}, \;\; \widetilde{\bs R} \in \R^{|S|\times N}.
\end{align*}
We observe that if we start with an $M$ vector of zeros, and insert into entries $S$ any nontrivial $|S|$-vector from the range of $\widetilde{\bs{Q}}_2$, then this constructed vector is in the cokernel of $\bs{V}$. Hence, we have constructed a cokernel vector of $\bs{V}$ requiring only $\mathcal{O}(|S| N)$ storage. If $|S| = N+1$, this reduces the storage requirement to $\mathcal{O}(N^2)$. In this streaming setting, we do form $\widetilde{\bs{Q}}_2$ explicitly, but it is small: its dimensions are $|S| \times (|S|-N)$, which is $(N+1) \times 1$ when $|S| = N+1$ (we focus on this in what follows). Forming the QR decomposition of $\bs{V}_{S\ast}$ costs $\mathcal{O}(|S| N^2)$ when $|S| = \mathcal{O}(N)$, which is the regime we consider.

A straightforward extension of the above idea to an iterative version of a \texttt{CSP} algorithm would order the elements in $[M]$ in any way, say encoded as a vector $\Sigma \in [M]^M$ whose values are a permutation of the elements of $[M]$, and for some fixed $k$ independent of $M$ and $N$ (chosen small for reduced runtime and storage complexity; we later focus on $k=1$), initialize $S$ as the first $N+k$ elements of this ordering and then repeatedly prune one index, and then add another. Since at each step at least one node is removed, the procedure still terminates in at most $M-N$ steps. We denote this \textit{streaming} variant of \Caratheodory-Steinitz pruning the \texttt{SCSP} algorithm, shown in \Cref{alg:scsp}. That the rows of $\bs{V}$ can be streamed to accelerate computation was also observed in \cite{vandenbos_generating_2020}. In addition, a streaming variant of measure compression specialized to a finite-size low discrepancy sequence was proposed in \cite{Elefante2022CQMCAI}.

\begin{algorithm}
  \caption{\texttt{SCSP}: Streaming \Caratheodory-Steinitz pruning}\label{alg:scsp}
    {\raggedright\noindent
      \textbf{Input: } $\bs{V} \in \mathbb{R}^{M\times N}$, $\bs{w} \in \mathbb{R}^M_{++}$, $k \in \N$, $\Sigma \in [M]^M$\\
        \textbf{Output: } $S$ with $|S| \leq N$, $\bs{w}_S \in \R^{|S|}_{++}$\par}
    \begin{algorithmic}[1]
        \State $S = \Sigma_{[N+k]}$, pop first $N+k$ indices from $\Sigma$. 
        \While{$|S| > N$}
          \State Compute ${\bs{n}} \in \mathrm{ker}({\bs{V}_{S*}}^T)$ \Comment{$\mathcal{O}((N+k)N^2)$}\label{lst:scsp-kernel-vector}
          \State Identify $S_\pm$ and compute $c_{\pm}$ in \eqref{eq:mc-def} using $S_{\pm}$, $\bs{w}_S$, $\bs{n}$.
          \State Choose $\sigma \in \{+,-\}$ as in \eqref{eq:sigma} \Comment{\texttt{SigSelect}, e.g., as in \eqref{eq:sigselect-example}}
          \State Set $\bs{w}_S \gets \bs{w}_S - c_\sigma \bs{n}$, let $P = \left\{ q \in S \;\;\big|\;\; w_q = 0 \right\}$. \label{lst:scsp-weight-prune}
          \State $S \gets S\backslash P$, pop first $\min(|P|,|\Sigma|)$ elements of $\Sigma$ and add to $S$. \label{lst:scsp-augment}
        \EndWhile
    \end{algorithmic}
\end{algorithm}

The \scsp{} algorithm now requires only $\mathcal{O}((N+k)N)$ storage, since only $N+k < M$ rows of the full matrix $\bs{V}$ and full initial weight vector $\bs{w}$ are stored at a time. The computational complexity of the \scsp{} algorithm is $\mathcal{O}((M-N)(N+k) N^2)$ because we expend $\mathcal{O}((N+k)N^2)$ effort to compute a cokernel vector of $\bs{V}$ at most $M-N$ times. Moving forward, we will assume fixed, small $k$ in which case the storage complexity is $\mathcal{O}(N^2)$ and the computational complexity is $\mathcal{O}((M-N)N^3)$. We reduce the complexity's polynomial order on $N$ by 1 in the next section.

\subsection{The \gscsp{} algorithm: \texorpdfstring{$\mathcal{O}(N^2)$}{O(N\^{2})} per-iteration complexity}\label{ssec:gscsp}
The computational bottleneck in a straightforward implementation of \Cref{alg:scsp} is the $\mathcal{O}(N^3)$ complexity of line \ref{lst:scsp-kernel-vector} that computes a cokernel vector of $\bs{V}$. At the first iteration, this cost is necessary, but at subsequent iterations the current iteration's matrix $\bs{V}$ differs from the previous iteration's by $|P|$ rows (typically one row); we can therefore employ low-rank modifications of the previous iterate's QR decomposition to generate the QR decomposition of the current iterate. More precisely, we first \textit{downdate} the previous iterate's QR decomposition by removing $|P|$ rows from $\bs{V}$, and then update the QR decomposition by adding back at most $|P|$ rows to $\bs{V}$. Efficient $\mathcal{O}(N^2)$ implementations of these procedures through Givens rotations are described in \cite{golub_matrix_1996}; we summarize the procedure here assuming $|P|=1$.

Consider $\bs{V} \in \R^{(N+k) \times N}$ corresponding to the previous iterate, and that $\bs{V}$ has the full QR decomposition $\bs{V} = \bs{Q} \bs{R}$. We let $\bs{G}$ denote a generic Givens rotation of size $N+k$.
To efficiently downdate the QR decomposition by removing row $i_{\text{rem}}$, one seeks to transform row $i_{\text{rem}}$ of $\bs{Q}$ to the vector $\pm \bs{e}_{i_{\text{rem}}}$ by building $N+k-1$ Givens rotations that zero out elements of this row:
\begin{align}\nonumber
  \bs{V} = \bs{Q}\bs{R} &= (\bs{Q}\bs{G}_{N+k-1}\cdots\bs{G}_{1})({\bs{G}_{1}}^T\cdots{\bs{G}_{N+k-1}}^T\bs{R})\\\label{eq:qr-downdate}
    &= \begin{bmatrix}
        & \vert & \\
        \bs{Q}_1 & \bs{0} & \bs{Q}_2\\
        & \vert & \\
        \text{-- }\bs{0}^T\text{ --} & \pm1 & \text{-- }\bs{0}^T\text{ --}\\
        & \vert & \\
        \bs{Q}_3 & \bs{0} & \bs{Q}_4\\
        & \vert & \\
    \end{bmatrix}\begin{bmatrix}
        \\
        \bs{R}_1\\
        \\
        \text{-- }\pm\bs{v}^T\text{ --}\\
        \\
        \bs{R}_2\\
        \phantom{0} % To make sure last row counts
    \end{bmatrix}
    = \begin{bmatrix}
        \bs{Q}_1\bs{R}_1 + \bs{Q}_2\bs{R}_2\\
        \text{--- }\bs{v}^T\text{ ---}\\
        \bs{Q}_3\bs{R}_1 + \bs{Q}_4\bs{R}_2\\
    \end{bmatrix},
\end{align}
where column $i_{\text{rem}}$ also equals $\pm \bs{e}_{i_{\text{rem}}}$ because the product of unitary matrices ($\bs{Q}$ and Givens rotations) is also unitary. The vector $\bs{v}^T$ coincides with row $i_{\text{rem}}$ of $\bs{V}$, $\bs{V}_{\{i_\text{rem}\}*}$. Then, letting $T = [N+k]\backslash \{i_{\text{rem}}\}$, by removing row $i_{\text{rem}}$ from the above expression, we have,
\[\bs{V}_{T*} = \widetilde{\bs{Q}}\widetilde{\bs{R}} = \begin{bmatrix}
    \bs{Q}_1 & \bs{Q}_2\\
    \bs{Q}_3 & \bs{Q}_4
\end{bmatrix}\begin{bmatrix}
    \bs{R}_1\\
    \bs{R}_2
\end{bmatrix} = \begin{bmatrix}
    \bs{Q}_1\bs{R}_1 + \bs{Q}_2\bs{R}_2\\
    \bs{Q}_3\bs{R}_1 + \bs{Q}_4\bs{R}_2\\
\end{bmatrix},\]
where $\widetilde{\bs{Q}}$ is unitary and $\widetilde{\bs{R}}$ remains upper triangular through proper ordering of the Givens rotations. Hence, \eqref{eq:qr-downdate} uses $\mathcal{O}(N^2)$ complexity to remove row $i_{\text{rem}}$ from a QR decomposition. See Section 3 of the Supplementary Information for pseudocode.

The second step is now to replace $\bs{v}$ in \eqref{eq:qr-downdate} with a new row vector, say $\widetilde{\bs{v}}^T$. Note that \eqref{eq:qr-downdate} with $\bs{v}$ replaced with $\widetilde{\bs{v}}$ is in a $\bs{Q} \bs{R}$ product form, but $\bs{R}$ is not upper triangular because the row $i_\text{rem}$ is dense. Again, we exercise Givens rotations to rectify this, first by zeroing out the first $i_{\textrm{rem}}-1$ entries of $\widetilde{\bs{v}}$, followed by ensuring the subdiagonal of $\bs{R}$ vanishes:
\begin{equation}\label{eq:qr-update}
 \begin{bmatrix}
    \bs{V}_1 \\
    \text{-- }\widetilde{\bs{v}}^T\text{ --}\\
    \bs{V}_2
\end{bmatrix} = \bs{Q}\bs{R} = \begin{bmatrix}
    & \vert & \\
    \bs{Q}_1 & \bs{0} & \bs{Q}_2\\
    & \vert & \\
    \text{-- }\bs{0}^T\text{ --} & \pm 1 & \text{-- }\bs{0}^T\text{ --}\\
    & \vert & \\
    \bs{Q}_3 & \bs{0} & \bs{Q}_4\\
    & \vert & \\
\end{bmatrix}{\bs{G}_{N}^T}\cdots{\bs{G}_1^T}\bs{G}_1\cdots\bs{G}_{N}\begin{bmatrix}
    \\
    \bs{R}_1\\
    \\
    \text{-- }\widetilde{\bs{v}}^T\text{ --}\\
    \\
    \bs{R}_2\\
    \phantom{0} % To make sure last row counts
\end{bmatrix}.
\end{equation}
The update procedure \eqref{eq:qr-update} requires $N$ Givens rotations for a complexity of $\mathcal{O}(N^2)$.
Therefore, this downdate-update procedure for fixed $k$ requires $\mathcal{O}(N^2)$ to compute a new kernel vector from the previous one. A formal pseudocode of the downdate-update procedure is presented in Supplementary Algorithm S1.

The \gscsp{} algorithm (``Givens \scsp'') is the augmentation of the \scsp{} algorithm by (i) retaining a dense $(N+k) \times N$ QR factorization throughout the process, (ii) using the downdate-update procedure described by \eqref{eq:qr-downdate} and \eqref{eq:qr-update} (implemented in Supplementary Algorithm S1) on line \ref{lst:scsp-augment} when updating the index set $S$, (iii) implementing line \ref{lst:scsp-kernel-vector} by simply slicing one of the trailing $k$ columns from the stored $\bs{Q}$ matrix.

The \gscsp{} algorithm is the proposed algorithm in this paper, which accomplishes \Caratheodory-Steinitz pruning with per-iteration $\mathcal{O}(N^2)$ complexity and storage. A summary of complexity and storage for all the algorithms discussed, assuming fixed $k$, is presented in \Cref{tab:asymptotics}. We again note that the NNLS-type procedure described in \cite{Elefante2022CQMCAI} may often outperform the expected NNLS runtime in \Cref{tab:asymptotics}, however, the total number of iterations for success is not known ahead of time.

The main computational advance of the \gscsp{} algorithm is efficient computation of a sequence of cokernel vectors for a matrix whose rows arrive in a streaming fashion. While we have exercised this advance to the particular situation of positive measure compression, this general-purpose procedure could have broader applications when similar numerical linear algebraic tasks are needed.

\begin{table}[h!]
    \centering
      \begin{tabular}{ l  c c }
         Algorithm &  Runtime & Storage  \\
         \hline
          \csp: \Cref{alg:csp}   & $(M-N) M N^2$ & $M N$ \\
          \scsp: \Cref{alg:scsp} & $(M-N) N^3$ & $N^2$ \\
         \gscsp$^{(\ast)}$: \Cref{alg:scsp}; Supplementary Algorithm S1 & $(M-N) N^2$  & $N^2$ \\
         \lp$^\dagger$ & $M^3$, $2^M$ & $MN$\\
         \nnls$^\dagger$ & $2^M$ & $MN$
        \end{tabular}
  \caption{$(M,N)$-asymptotic complexity of the algorithms presented in this paper. $M$ is the support size of $\mu_M$ and $N$ is the number of moments preserved with $M\geq N$. $(\ast)$: The first iteration of the \gscsp{} algorithm requires $\mathcal{O}(N^3)$ complexity to compute a dense QR factorization of an $(N+k) \times N$ matrix. ($\dagger)$: For \lp{} and \nnls{}, the approximate worst-case runtime complexities are presented, however, much faster convergence is observed in practice. For \lp{}, the two runtime complexities are for an interior point method and a simplex method, respectively.}
  \label{tab:asymptotics}
\end{table}

\section{Stability under measure perturbations}\label{sec:stability}
One numerical consideration is the effect of small perturbations of the original measure, $\mu_M$, on the resulting pruned quadrature rule. 
This can be important when considering issues of numerical round-off errors and in problems where the dense quadrature rule is an approximation or only approximate function evaluations are available.
Consider the polyhedron $W$ identified in \eqref{eq:W2}. Small perturbations of the weights $\bs{w}$ correspond to small perturbations of $W$. The addition of new nodes with small weights also continuously changes $W$ since this is equivalent to considering nodes with zero weights and perturbing them to slightly positive weights. The pruned quadrature rule is necessarily an extreme point of $W$ because it corresponds to at least $M-N$ active constraints (zero quadrature weights). Because continuous deformations of $W$ also continuously deform $\mathrm{ex}(W)$, the stability of the pruned quadrature rule with respect to small perturbations of the input quadrature rule is conceptually expected. Of course, the algorithmic way in which an element of $\mathrm{ex}(W)$ is identified may have different stability properties. 

We demonstrate explicitly in this section that the \scsp{} algorithm retains continuity of the pruned quadrature rule under small enough perturbations to the weights of the input quadrature rule or to the insertion of new nodes with small magnitude weights. This result immediately applies to the \gscsp{} algorithm since this algorithm is simply a computationally efficient version of \scsp. Instead of speaking in terms of quadrature rules, we will speak in terms of (discrete) measures. Consider the set of finite and finitely supported discrete signed measures on $X$, 
\begin{align}\label{eq:P-def}
  P &= \left\{ \nu\;\; \bigg|\;\; \nu = \sum_{m \in [M]} a_m \delta_{x_m}, \; M \in \N, \; a_m \in \R \textrm{  and  } x_m \in X \;\; \forall \;\; m\in[M]\right\}.
\end{align}
The set $P_+$ denotes the subset of $P$ that are non-negative measures:
\begin{align}\label{eq:Pplus-def}
  P_+ &= \left\{ \nu = \sum_{m \in [M]} a_m \delta_{x_m} \in P \;\;\bigg|\;\; a_m \geq 0\;\;\forall\;\; m \in [M]\right\}.
\end{align}
We compare two elements, $\alpha,\beta\in P_+$, using a variant of the total variation distance:
\begin{align}\label{eq:dTV}
  d_{\mathrm{TV}}(\alpha,\beta) &= \frac{|\alpha - \beta|}{|\alpha| + |\beta|},
\end{align}
where, assuming the nodes are distinct, $|\gamma|$ denotes the $\ell_1$ norm of the vector of weights of $\gamma \in P$:
\begin{align*}
  \gamma = \sum_{m \in [M]} d_m \delta_{x_m} \in P \hskip 10pt \Longrightarrow \hskip 10pt |\gamma| = \sum_{m \in [M]} |d_m|.
\end{align*}
If $\alpha$ and $\beta$ are both probability measures, then $d_{\mathrm{TV}}$ in \eqref{eq:dTV} reduces to the standard definition of total variation distance on discrete probability measures, which is $\frac{1}{2}$ times the $\ell_1$ norm difference between their mass functions. Over all probability measures having a fixed and finite maximum support size, $d_{\mathrm{TV}}$ is an equivalent distance to the Hellinger and Kullback-Leibler distance (when mass is bounded away from zero), and also to any translation-invariant distance induced by an $\ell^p$ norm over vectors. However, we reemphasize that this choice of norm is discontinuous with respect to the position of the nodes, so this choice of norm only captures perturbations to the weights and insertion of new nodes.

\subsection{Assumptions}
We require some assumptions on $\mu_M$ and the types of permissible measure perturbations. We first discuss a condition that allows us to ensure that the cokernel vectors used by the \scsp{} algorithm are well-behaved.
\begin{definition}
  The $N$-dimensional subspace $V$ is a \textit{Chebyshev system} with respect to $\mu_M$ if the Vandermonde-like matrix $\bs{V}$ in \eqref{eq:moment-conditions} satisfies $\det (\bs{V}_{S\ast}) \neq 0$ for every $S \subset [M]$ with $|S| = N$.
\end{definition}
For example, if $V$ is the subspace of degree-$(N-1)$ univariate polynomials, then it is always a Chebyshev system for any set $\mathrm{supp}(\mu_M)$ with at least $N$ distinct points. If $V$ is a subspace of multivariate polynomials and $\mu_M$ has its nodes drawn randomly with respect to some Lebesgue density function, then with probability one $V$ is a Chebyshev system with respect to $\mu_M$. The Chebyshev system property will enable us to assert uniqueness of cokernel vectors for submatrices of $\bs{V}$.

The \textit{ordering} of the nodes in the measure $\mu_M$ also plays a role in stability, and motivates the types of permissible perturbations to $\mu_M$. Recall that the \scsp{} algorithm starts from the $M \times N$ Vandermonde-like matrix $\bs{V}$ that is given as input. We first observe that the \scsp{} algorithm is \textit{not} stable with respect to permutations of the rows of $\bs{V}$. This can be conceptually understood by noting that permuting the rows of $\bs{V}$ would correspond to a change in the sequence of cokernel vectors $\bs{n}$ that is selected on line \ref{lst:scsp-kernel-vector} of \Cref{alg:scsp}. Hence, it is unreasonable to expect that the same pruned quadrature rule would result compared to the unpermuted case. From this observation, we recognize that it is not enough to discuss perturbations of a given measure $\mu_M$ under the total variation distance: we must also consider such perturbations subject to conditions that retain the ordering of the elements in $\mathrm{supp}(\mu_M)$.  In particular, we require enough assumptions so that the sequence of chosen cokernel vectors in the \scsp{} algorithm remains unchanged under perturbations. We will consider an ordering of the set $\mathrm{supp}(\mu_M) \subset X$; we denote this ordering as $\Sigma$: 
\begin{align*}
  \Sigma &\in \Pi(\mathrm{supp} (\mu_M)), & \Pi(Y) &= \mathrm{Sym}(Y) = \left\{\textrm{collection of permutations of $Y$}\right\}.
\end{align*}
Formally, we require the following assumptions on the measure $\mu_M$ and the hyperparameters of the \scsp{} algorithm.
\begin{assumption}\label{assum:scsp}
  Suppose the \scsp{} algorithm is run on $(V, \mu_M, \Sigma)$, for some $\Sigma \in \Pi(\mathrm{supp}(\mu_M))$. We assume that:
  \begin{itemize}
    \item $V$ is a Chebyshev system with respect to $\mu_M$.
    \item With $V$ fixed, then running the \scsp{} algorithm for any $(\mu_M,\Sigma)$ uses the same basis $v_n(\cdot)$ as in \eqref{eq:v-def}.
    \item $k = 1$, where $k$ is the integer input to \Cref{alg:scsp}.
    \item The \texttt{SigSelect} function is chosen as in \eqref{eq:sigselect-example}.
    \item The minimization problem \eqref{eq:sigselect-example} has a unique solution at every iteration.
  \end{itemize}
\end{assumption}
Given $(V,\mu_M,\Sigma)$ that satisfies Assumption \ref{assum:scsp}, note that running the \scsp{} algorithm generates a \textit{unique} output measure $\nu$ with \textit{unique} size-$N$ support $\mathrm{supp}(\nu)$. The uniqueness stems from the fact that Assumption \ref{assum:scsp} guarantees unique behavior of the algorithm:
\begin{itemize}
  \item Prescribing $\Sigma$ implies that the full Vandermonde-like matrix $\bs{V}$ is unique.
  \item $k=1$ implies that a cokernel vector from an $(N+1)\times N$ matrix, $\bs{V}_{S\ast}$ is computed at every step.
  \item $V$ being a Chebyshev system implies that $\mathrm{rank} (\bs{V}_{S\ast}) = N$ since any $N \times N$ submatrix must have full rank.
  \item The above two properties imply that the cokernel vector $\bs{n}$ at each iteration is unique up to multiplicative scaling. This scaling does not affect the algorithm: the updated weights $\bs{w} - c\,\bs{n}$, with $c$ chosen as in \eqref{eq:sigselect-example}, are unchanged if $\bs{n}$ is rescaled by a nonzero constant.
  \item Choosing \texttt{SigSelect} as in \eqref{eq:sigselect-example}, with the corresponding minimization problem having a unique solution, implies that the choice of pruned quadrature node is uniquely determined at every iteration.
\end{itemize}
Finally, we introduce the following set of admissible perturbations of $\mu_M$.
\begin{definition}\label{def:muM-perturbations}
  Let $(V, \mu_M, \Sigma)$ satisfy Assumption \ref{assum:scsp} and fix $\tau > 0$. 
  We define the set of admissible perturbations to $\mu_M$ as the set of measures and corresponding permutations on their supports, $(\widetilde{\mu}, \widetilde{\Sigma})$ for $\widetilde{\Sigma} \in \Pi(\mathrm{supp}(\widetilde{\mu}))$, as,
  \begin{align}\label{eq:muM-perturbations}
    P_\tau(\mu_M,\Sigma) \coloneqq \left\{ (\widetilde{\mu}, \widetilde{\Sigma}) \;\;\big|\;\; \widetilde{\mu} \in P_+, \;\; \widetilde{\Sigma}_{[M]} = \Sigma, \;\; \mathrm{supp}(\widetilde{\mu})\backslash \mathrm{supp}(\mu_M) \subset X_\tau \right\},
  \end{align}
  where $X_\tau = X_\tau(V)$ is defined as,
  \begin{align}\label{eq:Xtau}
    X_\tau \coloneqq \left\{x \in X \;\;\bigg|\;\; \sup_{v \in V\backslash\{0\}} \frac{|v(x)|}{\|v\|_{L_\mu^1(X)}} \leq \frac{1}{\tau} \right\}.
  \end{align}
\end{definition}
The set of valid measure perturbations therefore corresponds to measures that are positive, whose first ordered $M$ support points match the same ordered $M$ support points of the original measure, and whose support lies in the set $X_\tau \subseteq X$. The introduction of $\tau$ is a technical assumption; the precise value of $\tau$ is not conceptually important for theory, and it may be taken as an arbitrarily small, positive number. In particular, if the basis elements $v_j(\cdot)$ are all bounded over $X$, then there is a $\tau_0 > 0$ such that for any $\tau \in (0, \tau_0]$, we have $X_\tau = X$. Informally, $X_\tau$ exists to disallow nodal locations where $v(x)$ for arbitrary $v \in V$ has unbounded value.

\subsection{Stability results}
Our result on the stability of the \scsp{} algorithm is the following.
\begin{theorem}[\scsp{} and \gscsp{} stability]\label{thm:muM-perturbation}
  Let $(V, \mu_M,\Sigma)$ be given that satisfy Assumption \ref{assum:scsp}, and let $\nu = \scsp(\mu_M,V,\Sigma)$ be the output of the \scsp{} algorithm. For any fixed $\tau > 0$, define $P_\tau(\mu_M,\Sigma)$ as in \Cref{def:muM-perturbations}. Then the \scsp{} algorithm is locally Lipschitz (in particular continuous) with respect to the total variation distance in a $d_{\mathrm{TV}}$-neighborhood of $P_\tau(\mu_M,\Sigma)$ around $(\mu_M, \Sigma)$. I.e., there are positive constants $\delta_0 = \delta_0(\mu_M, \Sigma,\tau)$ and $C = C(\mu_M, \Sigma,\tau)$ such that for any $(\widetilde{\mu}_M, \widetilde{\Sigma}) \in P_\tau(\mu_M, \Sigma)$ satisfying $d_{\mathrm{TV}}(\mu_M, \widetilde{\mu}_M) < \delta_0$, then
  \begin{align*}
    d_{\mathrm{TV}}(\nu, \widetilde{\nu}) \leq C\, d_{\mathrm{TV}}(\mu_M, \widetilde{\mu}_M),
  \end{align*}
  where $\nu = \scsp(\mu_M,V,\Sigma)$ and $\widetilde{\nu} = \scsp(\widetilde{\mu}_M,V,\widetilde{\Sigma})$. 
\end{theorem}
See Appendix \ref{sec:scsp-proof} for the proof of this statement. The main message is that the \scsp{} (and hence also \gscsp) algorithm is robust to small enough perturbations (even perturbations that add nodes), provided those perturbations do not change the ordering of the original nodes.

A similar stability argument can be proven for the \nnls{} algorithm; the formal assumption, theorem, and proof are given in the Supplementary Information. The main difference is that \nnls{} is agnostic to ordering of the original nodes, but our analysis does not allow the addition of new nodes. In particular, fixing $\mu_M \in P_+$, we define the following set of perturbations formed by simply modifying the existing weights:
\begin{align*}
  P_{\mathrm{NNLS}}(\mu_M) &\coloneqq \left\{ \mu \in P_+ \;\;\big|\;\; \mathrm{supp} (\mu) = \mathrm{supp} (\mu_M) \right\}.
\end{align*}
Under additional uniqueness and nondegeneracy assumptions, the \nnls{} algorithm is locally Lipschitz with respect to $d_{\mathrm{TV}}$ for perturbations in $P_{\mathrm{NNLS}}(\mu_M)$. However, the set of admissible \nnls{} perturbations $P_{\mathrm{NNLS}}$ is considerably more restrictive than the admissible \scsp{} perturbations $P_{\tau}$, since the latter allows adding new nodes whereas the former does not. Our numerical results indicate that \nnls{} is not robust to adding new nodes.

\section{Numerical Results}\label{sec:numerical}

We investigate the efficacy of the \gscsp{} algorithm in this section. All experiments are performed in Julia version 1.11.5. The runtimes are compared on a Dell PowerEdge R6515 server equipped with an AMD EPYC 7702P processor (2.0 GHz, 64 cores / 128 threads) and 256 GB of RAM hosted on 64-bit (x86\_64) Ubuntu 24.04 LTS. We compare the following algorithms:
\begin{itemize}
  \item[\gscsp] The proposed algorithm of this manuscript: \Cref{alg:scsp} with the efficient Givens up/downdating described in \Cref{ssec:gscsp} with $k=1$ and \texttt{SigSelect} as in \eqref{eq:sigselect-example}. Our implementation is in the package \texttt{CaratheodoryPruning.jl}.
  \item[\nnls] The non-negative least squares procedure described in \Cref{ssec:alternatives}. For implementation we use the package \texttt{NonNegLeastSquares.jl} (\texttt{alg=:nnls}), which makes use of a version of the Lawson-Hanson algorithm with Householder reflections to solve intermediate least squares problems \cite{lawson_1995}.
  \item[\lp] The linear programming approach described in \Cref{ssec:alternatives} and in \eqref{eq:lp}. Unless otherwise stated, the vector $\bs{c}$ in \eqref{eq:lp} is formed with uniform random entries between 0 and 1. The implementation we use is the package \texttt{JuMP.jl} \cite{lubin_2023} with the \texttt{HiGHS} solver.
\end{itemize}
Note that alternative NNLS-based algorithms exist, mainly aimed at accelerating computation \cite{bro1997fast,dessole_lawson_hanson}. We have compared these alternatives (using software accompanying the literature) and the package we report above is consistently among the fastest in our experience for a ``pure'' NNLS algorithm. Further accelerations are possible through meta-algorithms, e.g., \cite{Elefante2022CQMCAI} that exercises one of the procedures above as an ingredient, but such meta-algorithms are not logically comparable to the standalone procedures we identify above.

\subsection{Computational complexity}
\begin{figure}[htbp]
    \centering
    \includegraphics[height=0.255\textwidth]{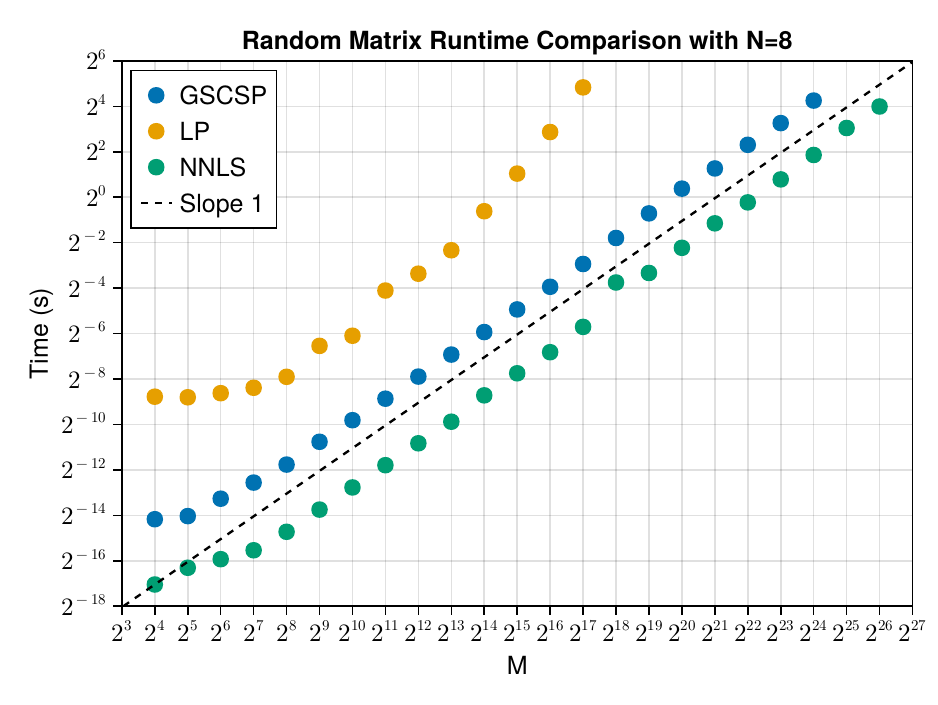}
    \includegraphics[height=0.255\textwidth,trim={1cm 0 0 0}, clip]{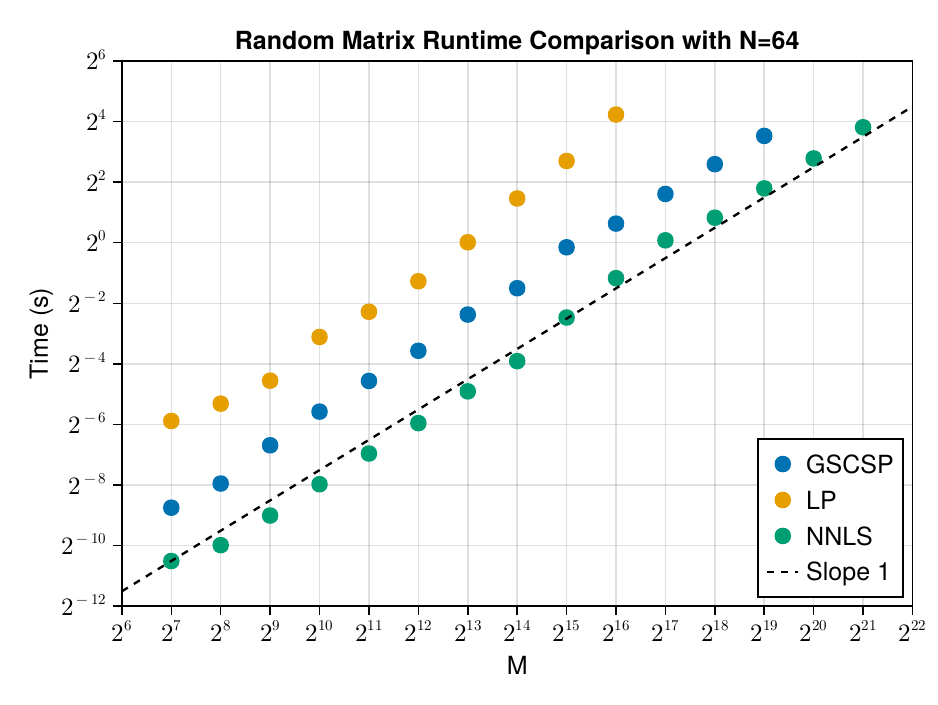}
    \includegraphics[height=0.255\textwidth,trim={1cm 0 0 0}, clip]{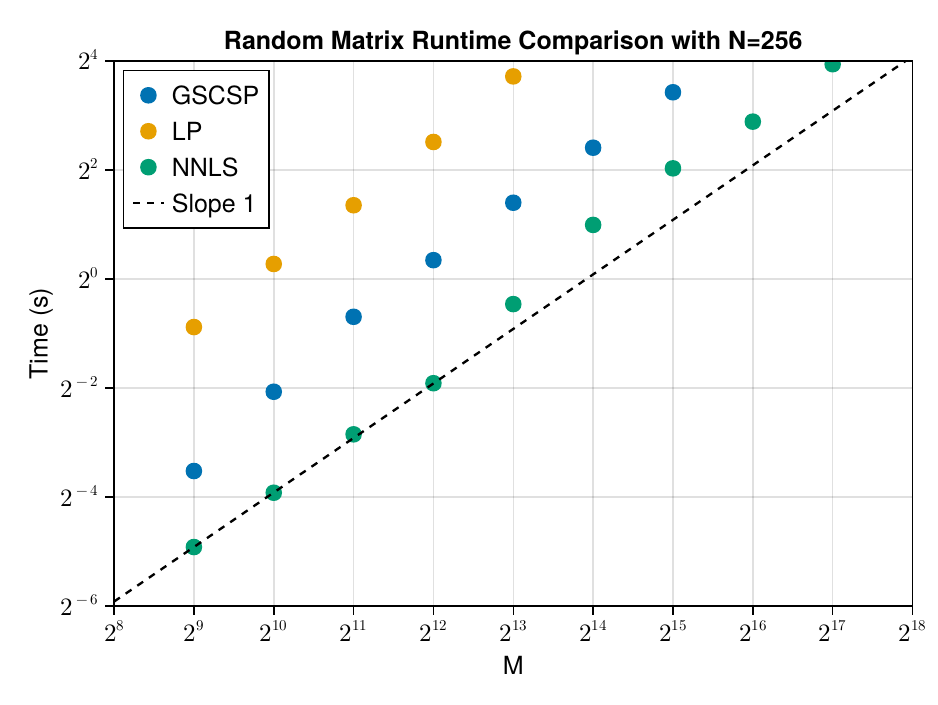}
    \caption{Runtime comparison of \gscsp{} pruning procedure to the \lp{} and \nnls{} approaches for $N=8$ (left), $N=64$ (middle), and $N=256$ (right) varying $M$. Each scatter point refers to a mean over 20 trials.}\label{fig:timings}
\end{figure}

We verify computational complexity in this section, in particular that \gscsp{} (or even just \scsp) requires linear complexity in $M$. For a given $(M,N)$, we generate $\bs{V}, \bs{w}$ as having independent and identically distributed (iid) entries uniformly between 0 and 1. \Cref{fig:timings} illustrates runtime comparisons of \gscsp, \lp, and \nnls{} approaches with fixed values of $N$ and increasing $M$. The results demonstrate that all methods have linear runtime complexity in $M$ and that in these regimes the \nnls{} procedure is the fastest followed by \gscsp{}, and then \lp. Increasing $N$ (\Cref{fig:timings}, middle and right) seems to have the effect of reducing the gap in runtime between \gscsp{} and \lp. Note that while \nnls{} is empirically the fastest algorithm in this experiment, the \gscsp{} algorithm we propose can easily operate in a streaming manner with predictable, asymptotically optimal runtime and requires less storage. We demonstrate this capability in the next section by taking very large $M$ that prevents dense storage of the matrix $\bs{V}$, which would be assumed impossible for straightforward implementations of \lp{} and \nnls{}.

\subsection{Quadrature on manufactured domains}
We present two sets of examples that demonstrate the flexibility of this approach in generating positive quadrature rules, and in particular in compressing very large quadrature rules. In all the examples of this section, we use the \gscsp{} algorithm to generate pruned quadrature. We will take $M = 10^9$. Let $X \subset \R^d$ be a $d$-dimensional compact domain; for visualization purposes we will focus on $d = 2, 3$. We consider the uniform probability measure $\mu$ on $X$, and we generate a random measure $\mu_M$ by sampling $M$ iid points $x_m$ from $\mu$ (by rejection sampling), and assign uniform weights $w_m = \frac{1}{M}$. Hence, while $\mu_M$ that we consider is random, we expect only small perturbations due to this randomness because the large $M = 10^9$ suggests we are well within the asymptotic regime of probabilistic concentration. In this large-$M$ regime, the alternative \lp{} and \nnls{} algorithms are simply infeasible to use due to computational storage requirements. With $M=10^9$ and $N\approx 70$, it would take $\approx$560GB of memory to store the dense Vandermonde matrix in double precision. Let $x = (x^{(1)}, \ldots, x^{(d)}) \in \R^d$ and $\alpha = (\alpha_1, \ldots, \alpha_d) \in \N_0^d$. We will use the following standard multi-index set definitions for $r\geq0$:
\begin{align*}
  A_{\mathrm{HC}, r} &\coloneqq \left\{ \alpha \in \N_0^d \;\;\big|\;\; \left\|\log (\alpha + 1) \right\|_1 \leq \log(r+1) \right\},\\
  A_{p, r} &\coloneqq \left\{ \alpha \in \N_0^d \;\;\big|\;\; \left\|\alpha\right\|_p \leq r\right\},\\
  A_{\mathrm{TD}, r} &\coloneqq A_{1,r}.
\end{align*}
In all our examples, a basis for $V$ is formed as a collection of $d$-fold products of univariate functions, where each basis function is constructed from one index $\alpha$ in a multi-index set. For example, each column of the Vandermonde-like matrix $\bs{V}$ is formed by choosing one $\alpha$, which defines a basis function $\phi \in V$ that is evaluated on $\{x_m\}_{m \in [M]}$.

\begin{table}[htbp]
\centering
\small
\begin{tabular}{@{}cllllc@{}}
\toprule
Example & Domain & $d$ & Basis family $\phi_{\alpha}(x)$ & Index set $A$ & Nodes \\
\midrule
$X_1$ & Mickey mouse shape & 2 & $H_{\alpha_1}(x^{(1)})H_{\alpha_2}(x^{(2)})$ & $A_{\mathrm{HC},20}$ & 70 \\
$X_2$ & Pumpkin shape & 2 & $J_{\alpha_1}(x^{(1)})J_{\alpha_2}(x^{(2)})$ & $A_{1/3,25}$ & 70 \\
$X_3$ & Spiral shape & 2 & $L_{\alpha_1}(x^{(1)})L_{\alpha_2}(x^{(2)})$ & $A_{\mathrm{TD},10}$ & 66 \\
$X_4$ & Torus volume & 3 & $\prod_{j=1}^3 (x^{(j)})^{\alpha_j}$ & $A_{\mathrm{HC},11}$ & 74 \\
\bottomrule
\end{tabular}
\caption{Manufactured-domain examples for $M=10^9$. Here $H_q$, $J_q$, and $L_q$ denote the univariate Hermite polynomial, Bessel function of the first kind, and Legendre polynomial, respectively.}
\label{tab:manufactured-domains}
\end{table}

\subsubsection{Two-dimensional domains, \texorpdfstring{$d=2$}{d=2}}\label{sssec:mc-d2}
\begin{figure}[htbp]
    \begin{subfigure}{0.305\textwidth}
    \includegraphics[trim={0 0 2.7cm 0},clip,width=\textwidth]{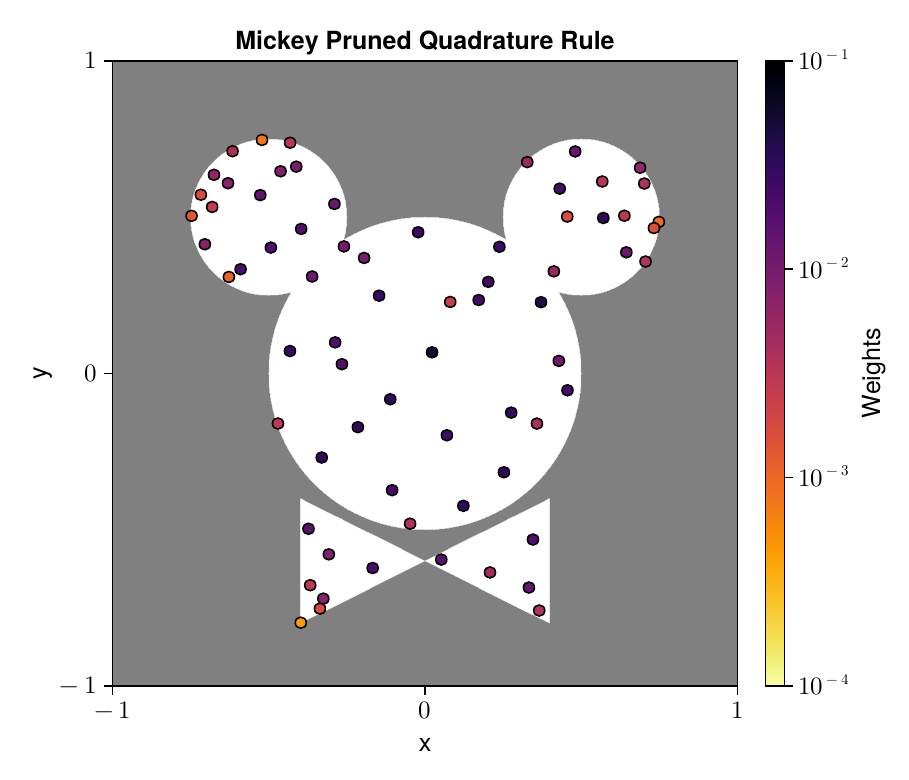}
    \end{subfigure}
    \begin{subfigure}{0.305\textwidth}
    \includegraphics[trim={0 0 2.7cm 0},clip,width=\textwidth]{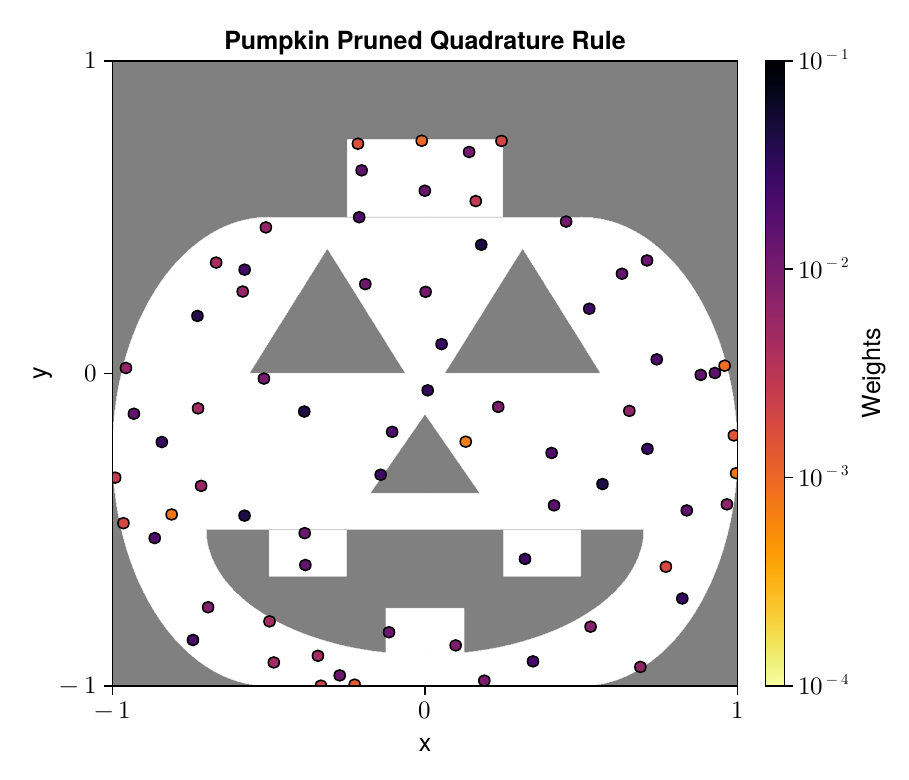}
    \end{subfigure}
    \begin{subfigure}{0.365\textwidth}
    \includegraphics[width=\textwidth]{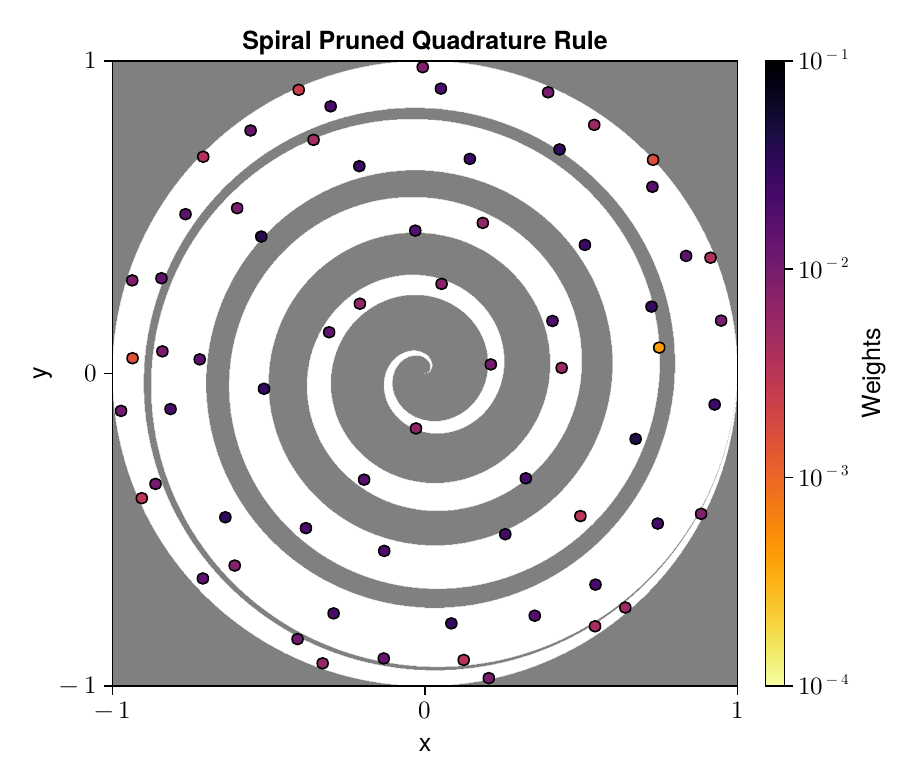}
    \end{subfigure}
    \caption{Examples of \gscsp{}-pruned quadrature rules on various 2D shapes. Left, middle, and right: examples $X_1$, $X_2$, and $X_3$, respectively.}
    \label{fig:2dpruning}
\end{figure}

\Cref{fig:2dpruning} illustrates the examples $X_1$-$X_3$ from \Cref{tab:manufactured-domains}: pruned Monte-Carlo integration rules on irregular shapes using product bases indexed by hyperbolic cross, quasi-norm, and total-degree sets. Some of these experiments are motivated by similar experiments in \cite{vandenbos_generating_2020}.
Note that none of the quadrature rules shown exhibit explicit symmetry that the subspace and/or domain possess. Even if we ensured that the full $M$-point rule were symmetric, pruned quadrature rules are generically not symmetric.

\subsubsection{Three-dimensional domains, \texorpdfstring{$d=3$}{d=3}}\label{sssec:mc-d3}
\begin{figure}[htbp]
    \centering
    \includegraphics[height=0.20\textheight]{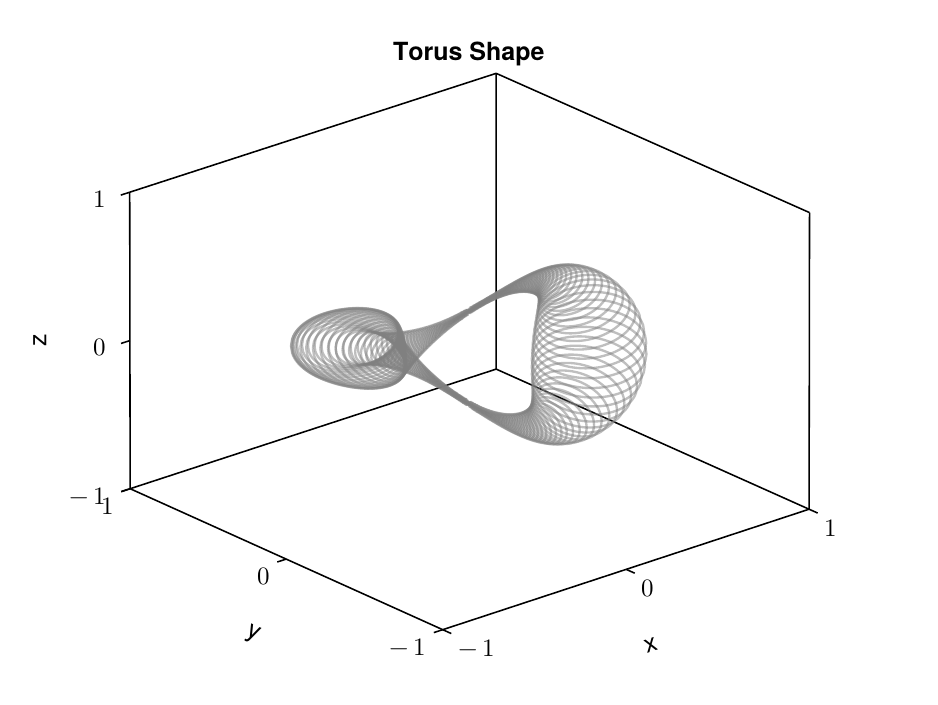}
    \includegraphics[height=0.20\textheight]{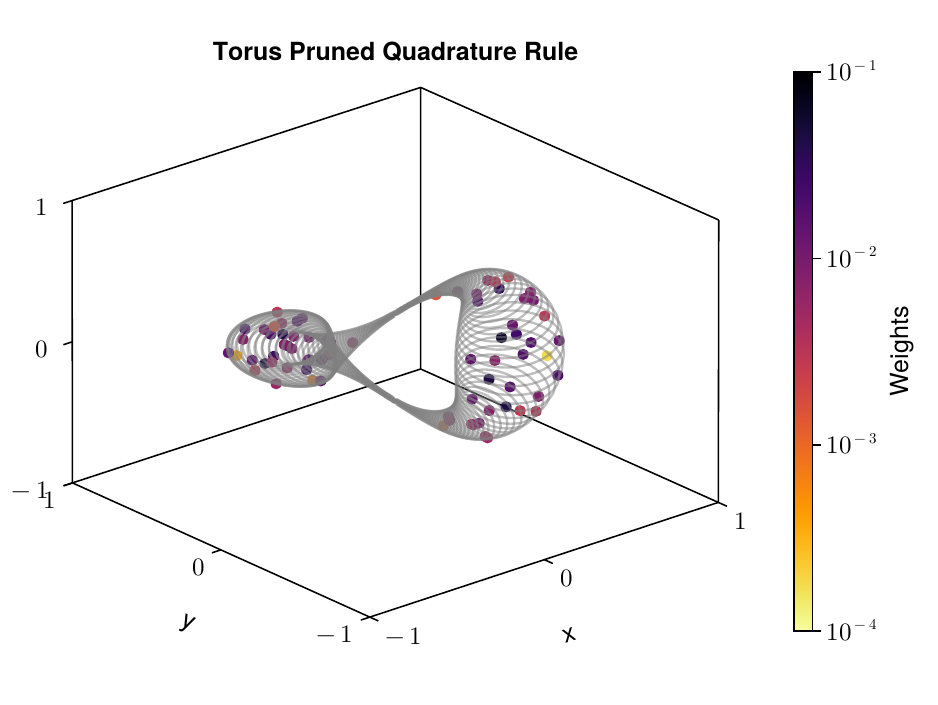}
    \caption{Example $X_4$ of a \gscsp{}-pruned quadrature rule on a 3D domain.}\label{fig:3dpruning}
\end{figure}

The dimension $d$ has no significant effect on the difficulty or complexity of the \gscsp{} algorithm; only the evaluation of the basis functions becomes slightly more expensive. \Cref{fig:3dpruning} illustrates $X_4$ from \Cref{tab:manufactured-domains}, the volume inside a torus whose radius and height vary with the polar angle in the two-dimensional plane.

\subsection{Stability of pruned quadrature rules}
We provide empirical experiments to complement the stability theory provided in \Cref{sec:stability}. In particular, we investigate $d_{\mathrm{TV}}(\nu_M, \widetilde{\nu}_M)$ when $\widetilde{\mu}_M$ is a small perturbation of $\mu_M$ and $\widetilde{\nu}_M$ is the result of pruning $\widetilde{\mu}_M$. We restrict ourselves to the model described in \Cref{sec:stability}, where $d_{\mathrm{TV}}(\mu_M, \widetilde{\mu}_M)$ is small subject to the constraint that ordering of nodes in $\widetilde{\mu}_M$ is fixed and that added nodes are appended to the existing ordering of the support of $\mu_M$. 

\begin{figure}
    \centering
    \includegraphics[trim={0 0 3cm 0},clip,height=0.28\textwidth]{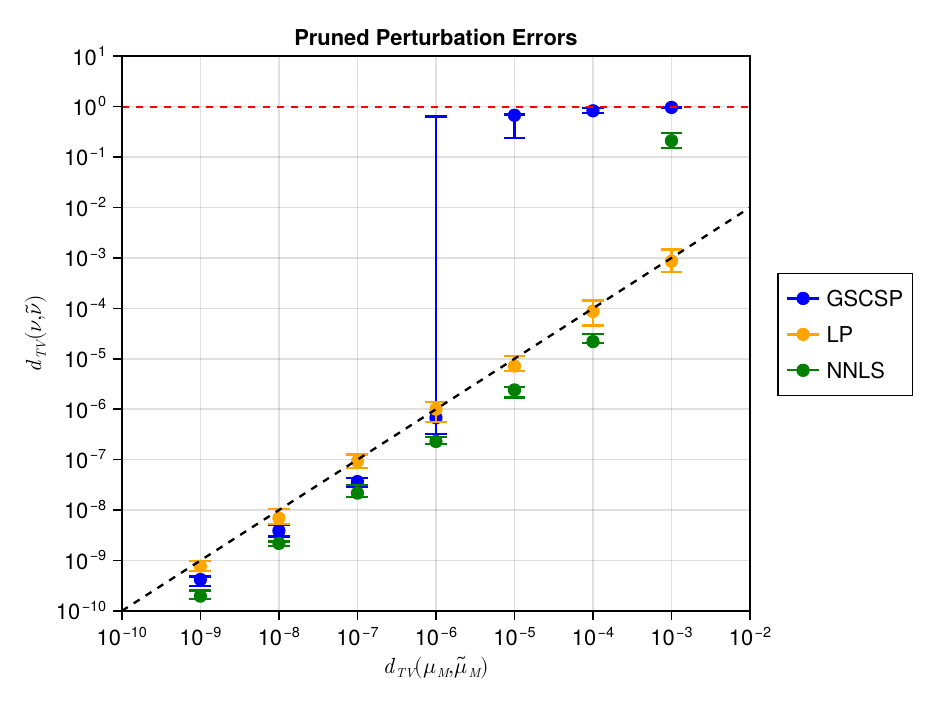}
    \includegraphics[trim={0 0 3cm 0},clip,height=0.28\textwidth]{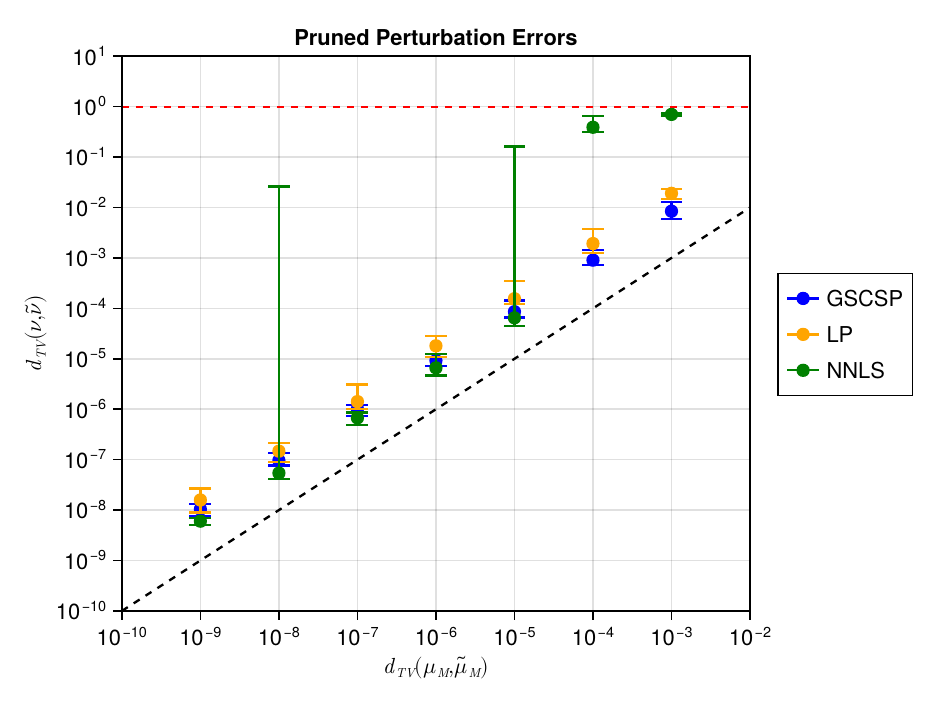}
    \includegraphics[height=0.28\textwidth]{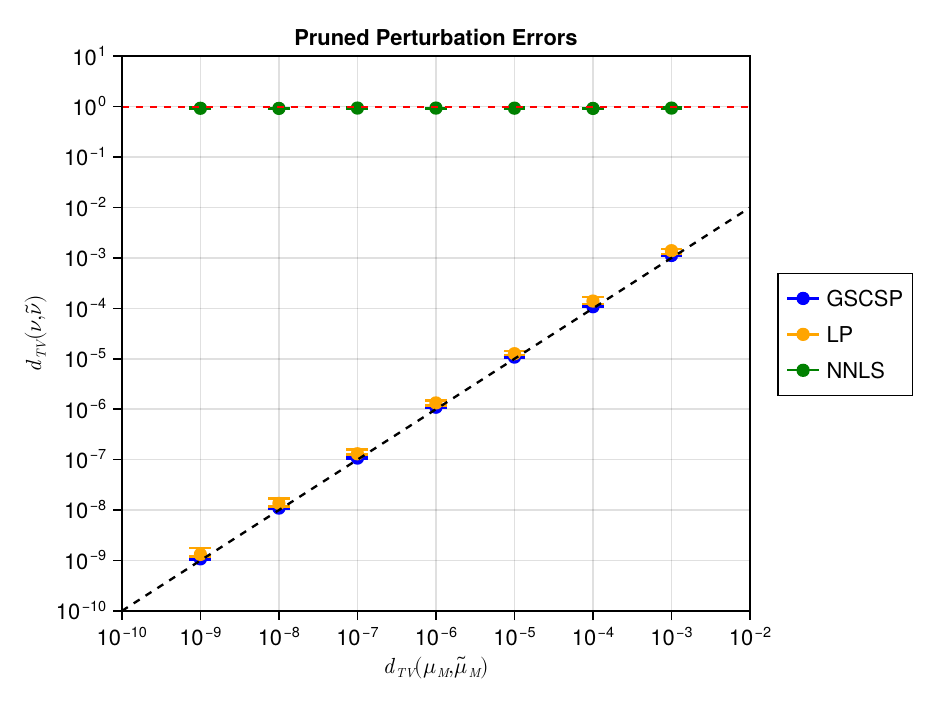}
    \caption{TV errors after various TV perturbations to a discrete measure $\mu_{M}$ with $M=10^4$. Left: perturbing the weights; middle: adding $10$ nodes with small positive weights; right: adding $10^4$ nodes with small positive weights.}
    \label{fig:prune_many}
\end{figure}

\Cref{fig:prune_many} illustrates the impacts of three different types of random perturbations on a fixed discrete measure, $\mu_M$ with $M=10^4$. In these examples, we form $\mu_M$ by rejection sampling points on the unit disk and assigning uniform initial weights $w_m = \frac{1}{M}$. The set $X$ and basis are chosen to be
\begin{align*}
    X &: \left\{x \in \R^2 \;\big|\; \|x\|\leq 1\right\} & V &= \mathrm{span}\left\{ L_{\alpha_1}(x^{(1)}) L_{\alpha_2}(x^{(2)}) \; \big|\; \alpha \in A_{\mathrm{HC},30} \right\},
\end{align*}
with $\dim(V) = 113$. Hence, each pruned rule has $N=113$ nodes.

We first compute $\nu_\mathrm{CSP}$, $\nu_\mathrm{LP}$, and $\nu_\mathrm{NNLS}$, by applying the \gscsp{}, \lp{}, and \nnls{} algorithms to $\mu_M$ respectively. Then, perturbations to achieve various total variation distances are applied to $\mu_M$ to form $\widetilde{\mu}_M$. We then compute $\widetilde{\nu}_\mathrm{CSP}$, $\widetilde{\nu}_\mathrm{LP}$, and $\widetilde{\nu}_\mathrm{NNLS}$, by applying the \gscsp{}, \lp{}, and \nnls{} algorithms to $\widetilde{\mu}_M$ respectively. Finally, we compute $d_\mathrm{TV}({\nu}_\mathrm{CSP},\widetilde{\nu}_\mathrm{CSP})$, $d_\mathrm{TV}({\nu}_\mathrm{LP},\widetilde{\nu}_\mathrm{LP})$, and $d_\mathrm{TV}({\nu}_\mathrm{NNLS},\widetilde{\nu}_\mathrm{NNLS})$. The scattered values are the median along with the 0.2 up to the 0.8 quantiles over 20 repetitions. To reduce randomness of the \lp{} algorithm, the vector $\bs{c}$ used in \eqref{eq:lp} is kept the same through all simulations with ones appended for appended nodes.

In the left panel of \Cref{fig:prune_many}, no new nodes are added to $\mu_M$, however, we apply random, mean zero, displacements to the weights of $\mu_M$ (maintaining positivity) to achieve various total variation perturbations. All methods are stable to this type of perturbation up to a certain magnitude displacement. (Note, however, that \gscsp{} has the smallest radius of stability in this example.) \lp{} was found to be the most stable in this case, followed by \nnls{}, followed by \gscsp{}. In the middle panel of \Cref{fig:prune_many}, the weights of $\mu_M$ are maintained and $10$ new (randomly sampled) nodes are inserted with uniform small weight to achieve various total-variation distances. In this case, the \nnls{} algorithm loses stability while the \gscsp{} algorithm seems to be slightly more stable than \lp{}. Finally, in the right panel of \Cref{fig:prune_many}, the same type of perturbation is used as in the middle panel, instead, $10^4$ nodes are appended. In this case, the \nnls{} algorithm is completely unstable while the \gscsp{} algorithm and \lp{} algorithms remain relatively stable. One reason for this may be that the magnitudes of the appended weights are inversely proportional to the number of appended nodes when considering the same total variation displacement between $\mu_M$ and $\widetilde\mu_M$.

\begin{figure}
    \centering
    \includegraphics[trim={0 0 3cm 0},clip,height=0.28\textwidth]{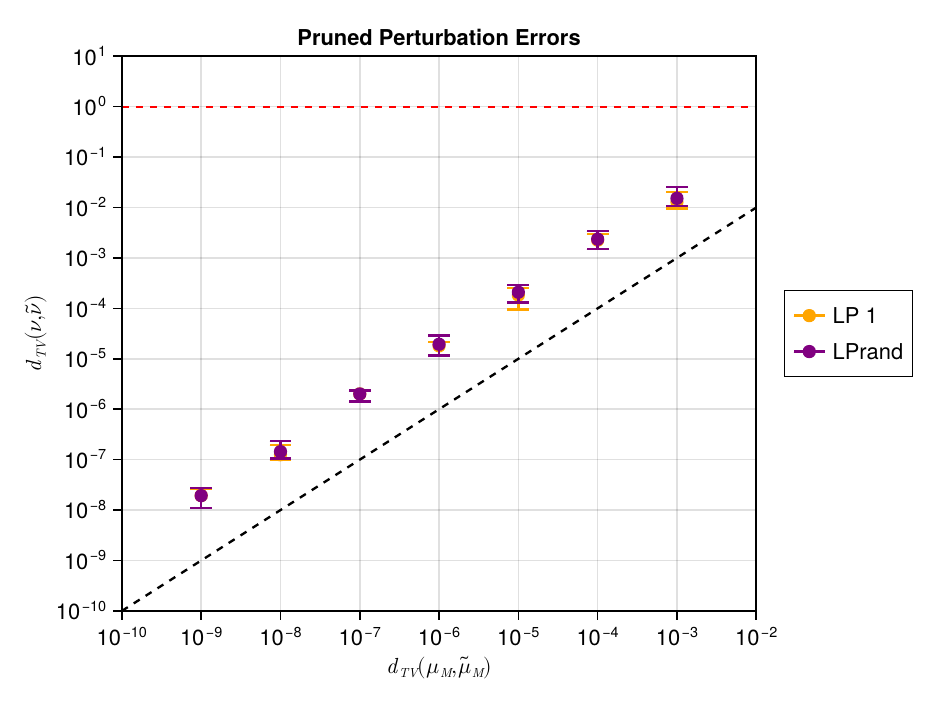}
    \includegraphics[trim={0 0 3cm 0},clip,height=0.28\textwidth]{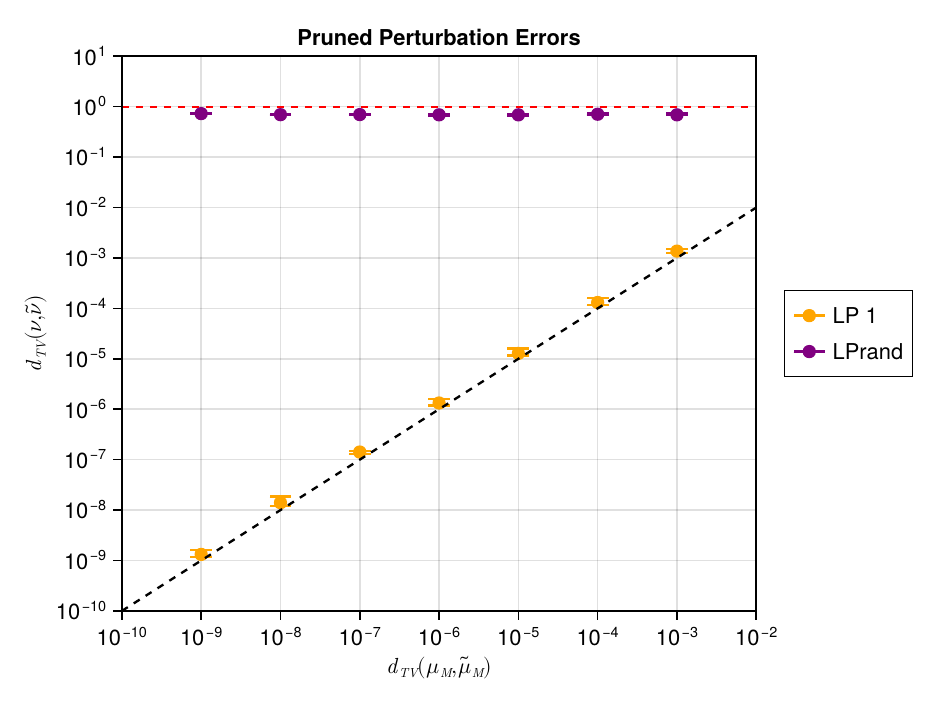}
    \includegraphics[trim={13cm 0 0 0},clip,height=0.28\textwidth]{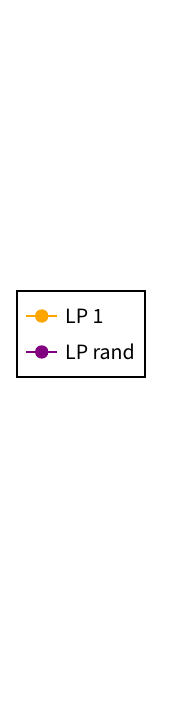}
    \caption{Replication of the middle and right panels of \Cref{fig:prune_many} with two choices of the \lp{} objective vector $\bs c$. LP 1: appending ones to $\bs c$; LP rand: appending $(0,1)$-uniformly sampled values to $\bs c$.}
    \label{fig:prune_many_lp}
\end{figure}

\Cref{fig:prune_many_lp} illustrates the same test as the middle and right panels of \Cref{fig:prune_many}, but compares the stable choice where ones are appended to $\bs{c}$ against the unstable choice where appended entries are sampled independently from the uniform distribution on $(0,1)$. In the left panel, where only 10 nodes are added, this choice did not have a significant impact on stability. However, the right panel of \Cref{fig:prune_many_lp} illustrates that in the situation of adding many new nodes, the choice of appending ones instead of random elements is crucial for stability.

\begin{figure}
    \begin{subfigure}{0.317\textwidth}
        \includegraphics[width=\textwidth]{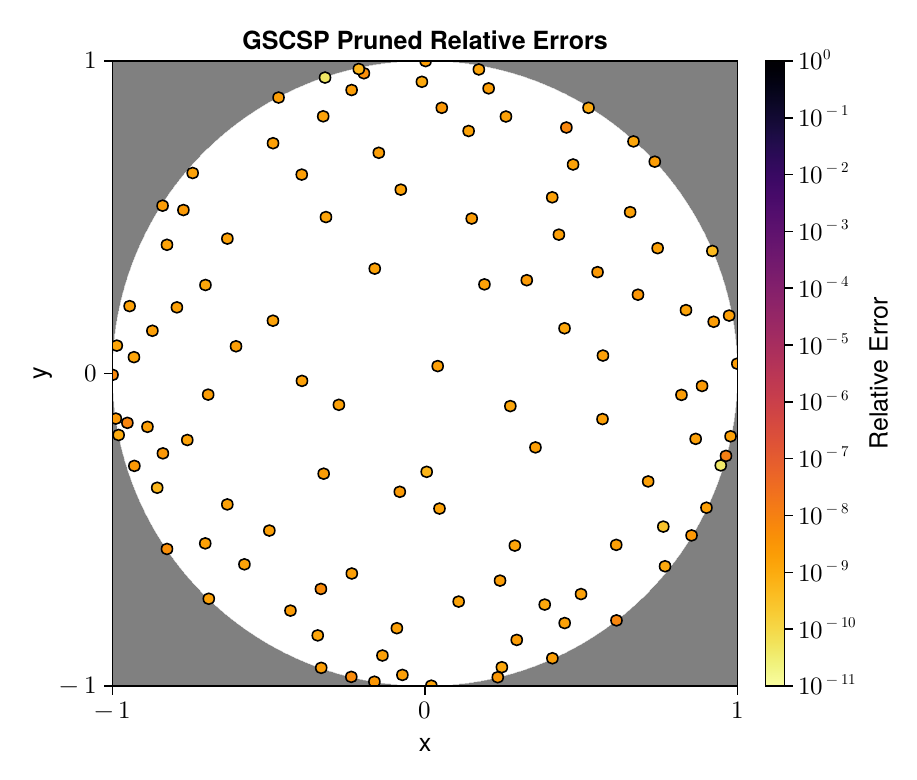}
        \includegraphics[trim={0 0 4cm 0},clip,width=0.83\textwidth]{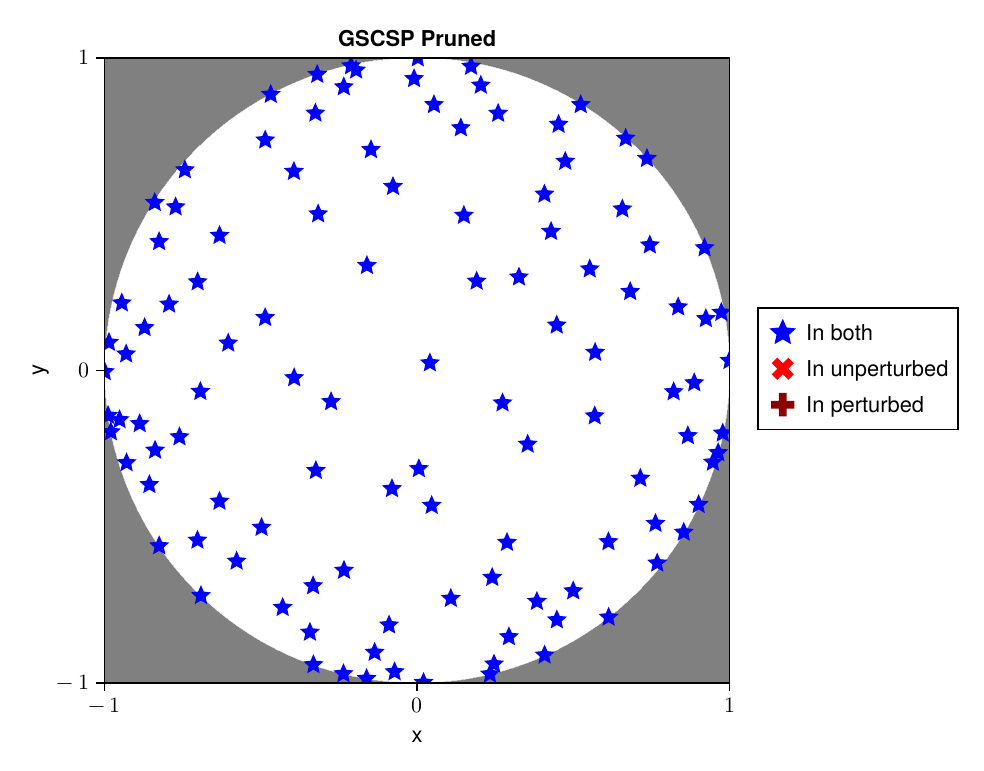}
    \end{subfigure}
    \begin{subfigure}{0.317\textwidth}
        \includegraphics[width=\textwidth]{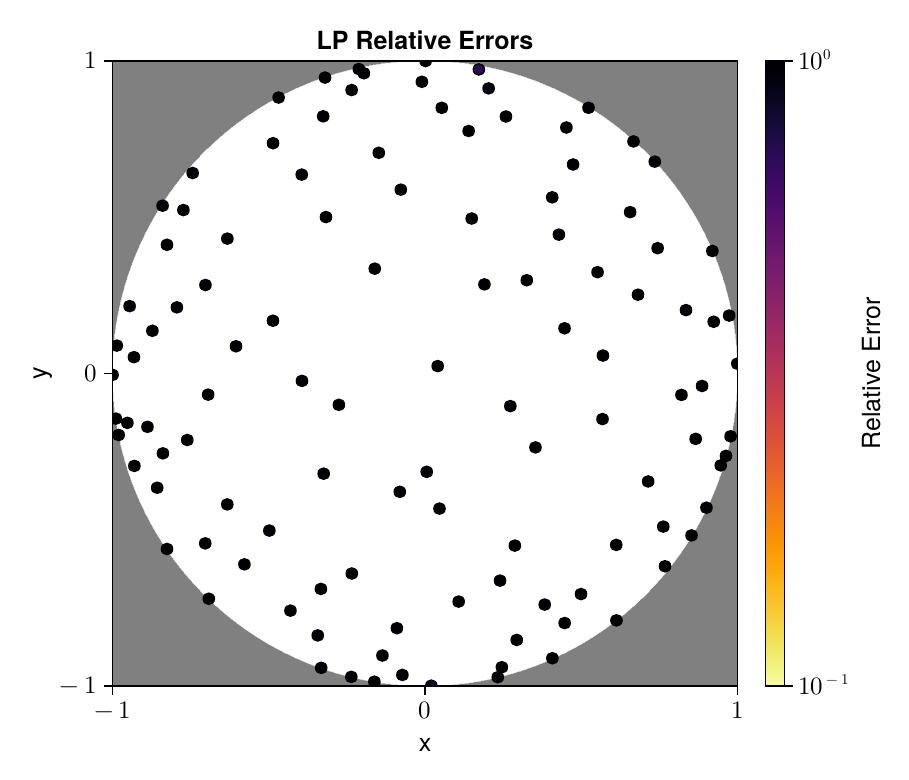}
        \includegraphics[trim={0 0 4cm 0},clip,width=0.83\textwidth]{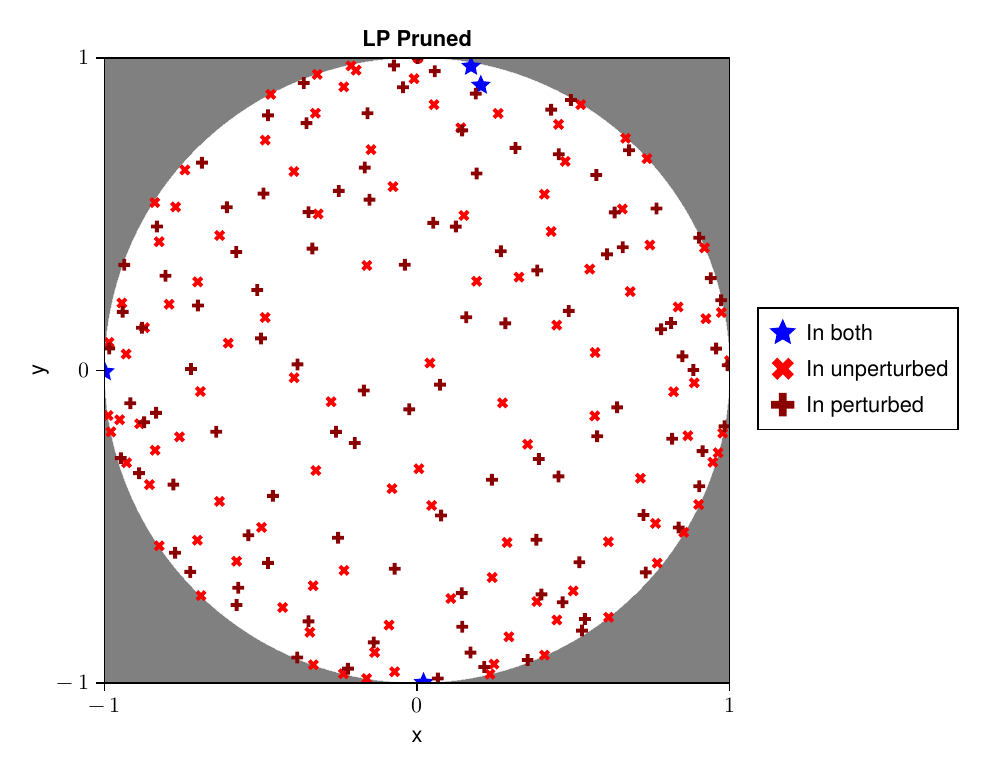}
    \end{subfigure}
    \begin{subfigure}{0.348\textwidth}
        \includegraphics[width=0.91092\textwidth]{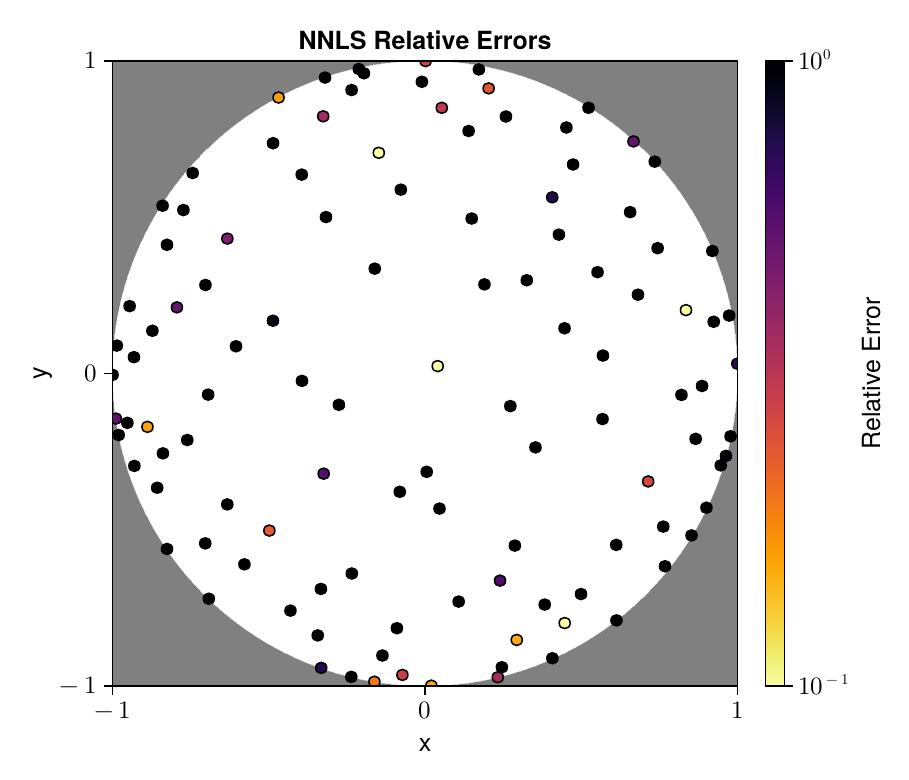}
        \includegraphics[width=\textwidth]{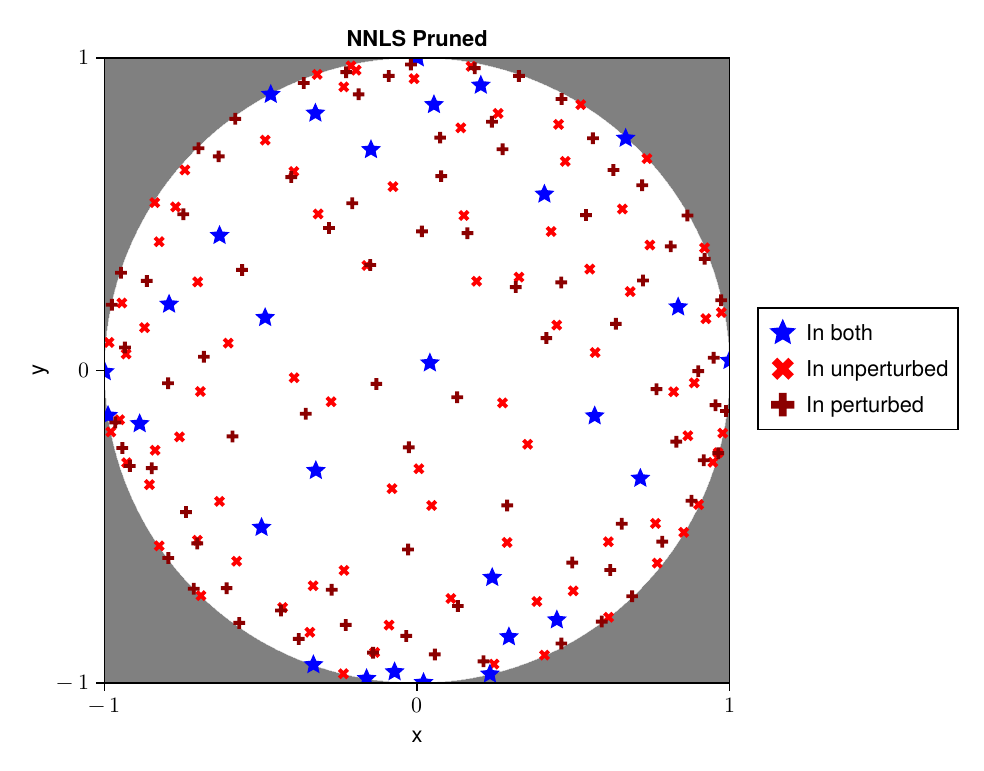}
    \end{subfigure}
    \caption{Relative errors and point comparison for pruned perturbed quadrature rules on unit disk. The reference pruned rule is generated using \nnls{} to prevent bias towards the \gscsp{} or \lp{} algorithms. The reference rule is then perturbed and re-pruned by each algorithm. The nodes of $\nu$ are displayed on the bottom row with blue $\star$ and red $\mathsf{x}$ characters; the nodes of $\widetilde{\nu}_\mathrm{CSP}$, $\widetilde{\nu}_\mathrm{LP}$, and $\widetilde{\nu}_\mathrm{NNLS}$ are displayed on the bottom row with blue $\star$ and maroon {$\scriptsize \boldsymbol{+}$} characters in each column respectively.}
    \label{fig:prunetest}
\end{figure}

In \Cref{fig:prunetest}, we visualize the instability of the \lp{} and \nnls{} methods. We use the same domain and basis as in the other stability test. In this test, we first form $\mu_M$ with $M=10^5$ points. We then prune the result down to $\nu$ using the \nnls{} algorithm. We know from former results that \nnls{} is not stable with respect to adding nodes, so the purpose of using the \nnls{} algorithm is to form a sparse quadrature rule which is not biased towards either the \gscsp{} or \lp{} algorithms. We then perturb the measure $\nu$ by adding $10^4$ nodes of uniform weight to form $\widetilde{\nu}$ such that $d_\mathrm{TV}(\nu,\widetilde{\nu}) = 10^{-9}$. Finally, we compute $\widetilde{\nu}_\mathrm{CSP}$, $\widetilde{\nu}_\mathrm{LP}$, and $\widetilde{\nu}_\mathrm{NNLS}$, by applying the \gscsp{}, \lp{}, and \nnls{} algorithms to $\widetilde{\nu}$ respectively. The top row of \Cref{fig:prunetest} displays the point-wise weight relative errors between the pruned rules and $\nu$. The bottom row of \Cref{fig:prunetest} illustrates which nodes are retained through this process and which are not. \Cref{fig:prunetest} clearly illustrates that the \gscsp{} algorithm is the most stable with respect to this type of perturbations.

\subsection{Application: cut-cell discontinuous Galerkin (DG) methods}

We conclude with an example of \Caratheodory-Steinitz pruning applied to the generation of quadrature rules for high order cut-cell discontinuous Galerkin (DG) methods. \Caratheodory-Steinitz pruning was previously used 
to construct reduced quadrature rules for cut-cell DG methods \cite{taylor2024entropy} and projection-based reduced order models (ROMs) \cite{qu2025entropy}. These reduced quadrature rules retain positivity while exactly satisfying certain moment conditions related to integration by parts \cite{chan2019skew}, and can be used to construct semi-discretely entropy stable discretizations for nonlinear conservation laws. We note that, for ROMs, a reference high cardinality quadrature rule is often available from the full order model. 

Here, we utilize reduced quadrature rules constructed using \Caratheodory-Steinitz pruning as described in \cite{taylor2024entropy} for high order cut-cell DG formulations of a 2D linear time-dependent advection-diffusion problem: 
\begin{equation}
\label{eq:advec_diff}
\frac{\partial u}{\partial t} + \nabla\cdot(\bs{\beta} u) - \epsilon\Delta u = f,\quad \bs{x} = (x,y) \in \Omega,
\end{equation}
where $\bs{\beta}(\bs{x}) = (-y, x)^T$ is a spatially varying advection vector, $\epsilon=10^{-2}$ is the diffusivity coefficient, $f(\bs{x}) = 1$ is a forcing term, and the domain, $\Omega$, is taken to be $[-1,1]^2 \setminus \Gamma$, where $\Gamma$ is the union of two circles of radius $R=0.4$ centered at $\frac{1}{2}(\mp1,\pm1)$.

On cut cells, the DG solution is often represented using physical frame total degree $P$ polynomials  \cite{giuliani2022two, taylor2025energy, taylor2024entropy}. Volume integrals over cut cells are typically challenging to compute due to the nature of the geometry and the physical-frame approximation space. In this work, we construct quadratures for volume integrals over cut cells which are \textit{exact} on the (isoparametrically approximated) cut cell for physical-frame polynomial integrands up to a certain degree.

To construct exact quadratures on cut cells, we first approximate cut elements using curved isoparametric subtriangulations.
Then, an exact quadrature on a cut cell can be constructed using a composite quadrature rule from simplicial quadratures of sufficient degree on each element of the curved subtriangulation. Because of the isoparametric representation of such a subtriangulation, one can show that a physical-frame polynomial integrand of degree $K$ can be exactly integrated using a quadrature rule of degree $K P + 2(P-1)$ on the reference element \cite{taylor2024entropy}.\footnote{In this example, surface integrals are exactly computed using standard one dimensional Gaussian quadrature rules of sufficient degree, but could also be pruned using a similar \Caratheodory-Steinitz procedure.} We emphasize that these quadrature rules are exact on the isoparametric approximation of the cut cell, not on the original cut cell, since the curved subtriangulation only approximates the cut-cell boundary. Pruning preserves this property: the pruned rule reproduces the moments of the unpruned rule exactly, and thus introduces no additional geometric error.

In this work, we take $K=2P-1$ such that the product of a degree $P$ polynomial and its derivative is exactly integrated. 
Thus, the reference quadrature rules we use to construct composite quadratures over cut cells should be exact for polynomials of degree $2P^2 + P-2$. These can be constructed using tensor products of one-dimensional Gaussian quadrature rules combined with a collapsed coordinate mapping \cite{karniadakis2013spectral}. For this problem, this construction requires $M=n_\mathrm{tri}(\lceil \frac{2P^2 + P-1}{2} \rceil+1)^2 \sim n_\mathrm{tri} P^4$ nodes in the dense, unpruned quadrature rule with $n_\mathrm{tri}=3$ subtriangles.
After mapping these to the physical domain, pruning is performed preserving the $N=P(2P+1) \sim 2P^2$ moments of the total degree set of degree $K=2P-1$ bivariate polynomials.
For a degree $P=7$ approximation, this implies the reference quadrature rule must be exact for an extraordinarily high polynomial degree of $103$, which results in a reference quadrature rule with $M=8427$ nodes. After pruning, we are left with a size $N=105$ point quadrature rule. Since interpolation to quadrature nodes and $L^2$ projection back to the basis are dense matrix-vector products whose cost is linear in the number of nodes, pruning reduces the cost of these operations by a factor of $M/N \approx 80$ on every cut element at each time-stepping stage.
\begin{figure}
\centering
\includegraphics[width=0.4\textwidth, trim={10em 2em 10em 2em},clip]{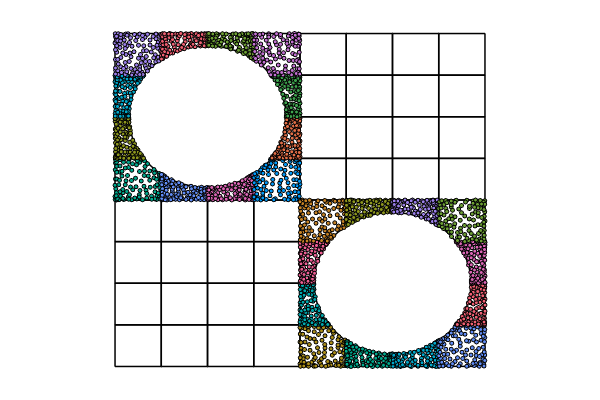}
\hspace{2em}
\includegraphics[width=0.4\textwidth, trim={10em 2em 10em 2em},clip]{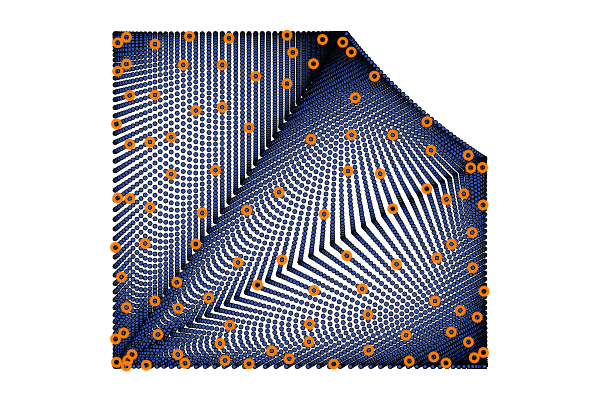}
\caption{The cut domain for the advection-diffusion problem (left) and pruned quadrature nodes compared with unpruned quadrature nodes on a single cut cell (right).}
\label{fig:cut_domain_nodes}
\end{figure}
\Cref{fig:cut_domain_nodes} shows the domain with pruned cut-cell quadrature points overlaid, as well as a comparison of the pruned and original unpruned quadrature rule with $P=7$. We note that, for the linear advection-diffusion problem in this paper, a positive moment-preserving quadrature rule is not strictly necessary to guarantee stability \cite{giuliani2022two, taylor2025energy}. However, the use of \Caratheodory-Steinitz pruning does ensure positive-definiteness of the mass matrix (via positivity) and high order accuracy (via exact satisfaction of moment conditions).

We first verify the high order accuracy of the cut-cell DG discretization. We solve \eqref{eq:advec_diff} with a constant advection $\boldsymbol{\beta} = (1.5, 1)^T$, diffusivity $\epsilon = 10^{-2}$, and forcing $f = 0$. We use a plane wave exact solution
\begin{equation}
u(\boldsymbol{x},t) = 1 + \frac{1}{2}\sin\big(\pi(x + y - (\beta_1 + \beta_2)t)\big)\,e^{-2\epsilon\pi^2 t},
\label{eq:cutcell_manufactured}
\end{equation}
from which boundary data is also generated. For the advective discretization, the exact solution is imposed at inflow boundaries where $\boldsymbol{\beta}\cdot\boldsymbol{n} < 0$, while for the viscous discretization, the exact solution is imposed as Dirichlet data on all of $\partial\Omega$. We compute solutions for degrees $P= 1, \ldots, 4$ on background Cartesian grids with mesh sizes $h = 1/2, 1/4, 1/8, 1/16$. The semi-discrete system is integrated to final time $T = 1$ using the adaptive 5(4) explicit Runge--Kutta method of Tsitouras with absolute and relative tolerances of $10^{-8}$. Table~\ref{tab:cutcell_convergence} reports $L^2$ errors at $T = 1$, which are computed using the volume quadrature rule on each cell. The average rates are close to the optimal rate of $O(h^{P+1})$.

\begin{table}[ht]
\centering
\footnotesize
\setlength{\tabcolsep}{3pt}
\begin{tabular}{c|cc|cc|cc|cc}
\hline
 & \multicolumn{2}{c|}{$P = 1$} & \multicolumn{2}{c|}{$P = 2$} & \multicolumn{2}{c|}{$P = 3$} & \multicolumn{2}{c}{$P = 4$} \\
$h$ & $L^2$ error & Rate & $L^2$ error & Rate & $L^2$ error & Rate & $L^2$ error & Rate \\
\hline
$1/2$ & $1.64\times 10^{-1}$ & -- & $1.25\times 10^{-2}$ & -- & $2.70\times 10^{-3}$ & -- & $1.16\times 10^{-4}$ & -- \\
$1/4$ & $3.64\times 10^{-2}$ & 2.17 & $1.95\times 10^{-3}$ & 2.68 & $3.26\times 10^{-4}$ & 3.05 & $6.07\times 10^{-6}$ & 4.25 \\
$1/8$ & $7.79\times 10^{-3}$ & 2.23 & $2.88\times 10^{-4}$ & 2.76 & $9.68\times 10^{-6}$ & 5.07 & $3.38\times 10^{-7}$ & 4.17 \\
$1/16$ & $1.73\times 10^{-3}$ & 2.17 & $2.21\times 10^{-5}$ & 3.70 & $5.73\times 10^{-7}$ & 4.08 & $6.28\times 10^{-9}$ & 5.75 \\
\hline
\end{tabular}
\caption{$L^2$ errors and observed convergence rates at $T = 1$ for the degree $P$ cut-cell DG discretization of the advection-diffusion problem with exact solution \eqref{eq:cutcell_manufactured}.}
\label{tab:cutcell_convergence}
\end{table}

\begin{figure}
\centering
\includegraphics[width=0.3\textwidth, trim={60em 2em 60em 2em},clip]{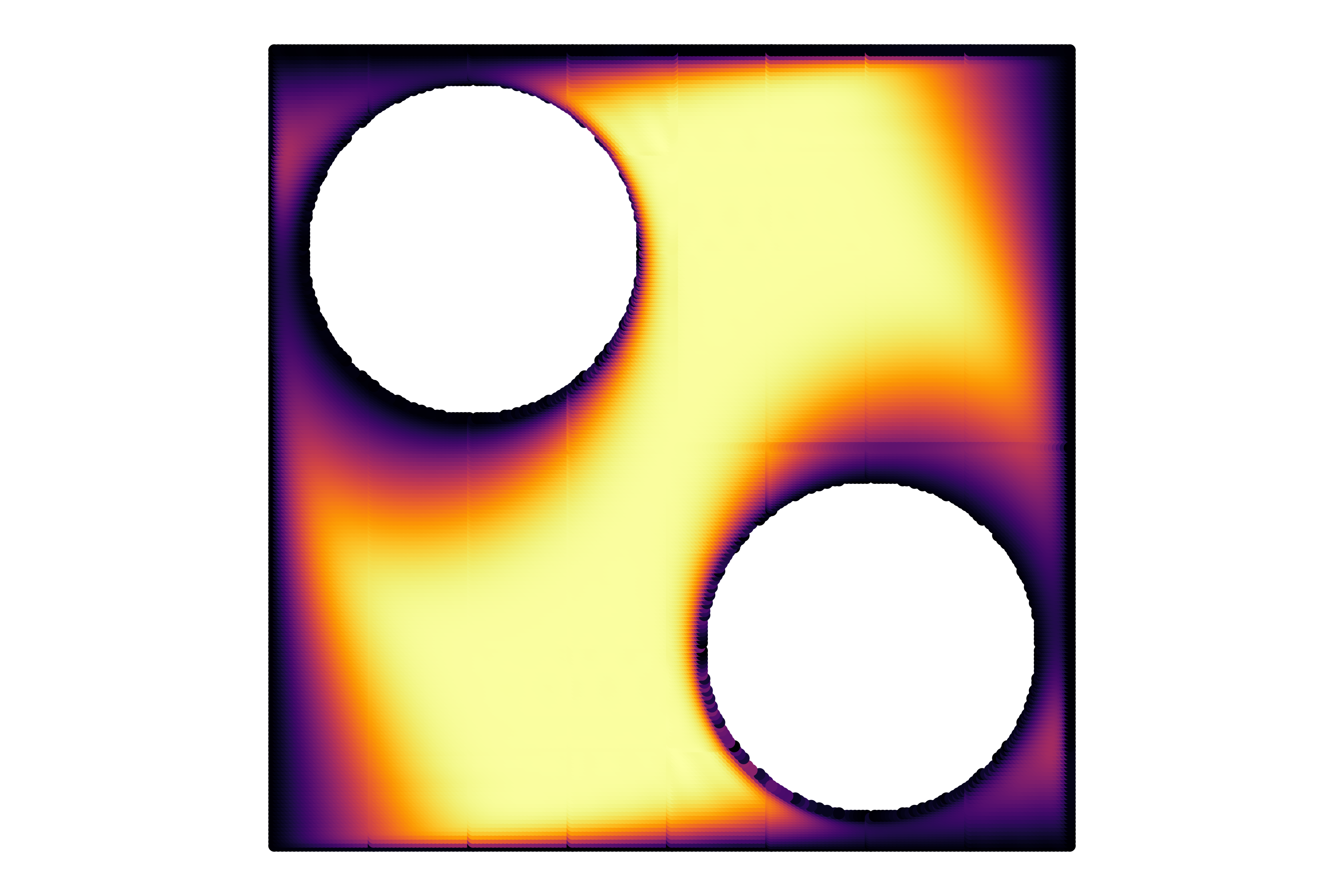}
\includegraphics[width=0.3\textwidth, trim={60em 2em 60em 2em},clip]{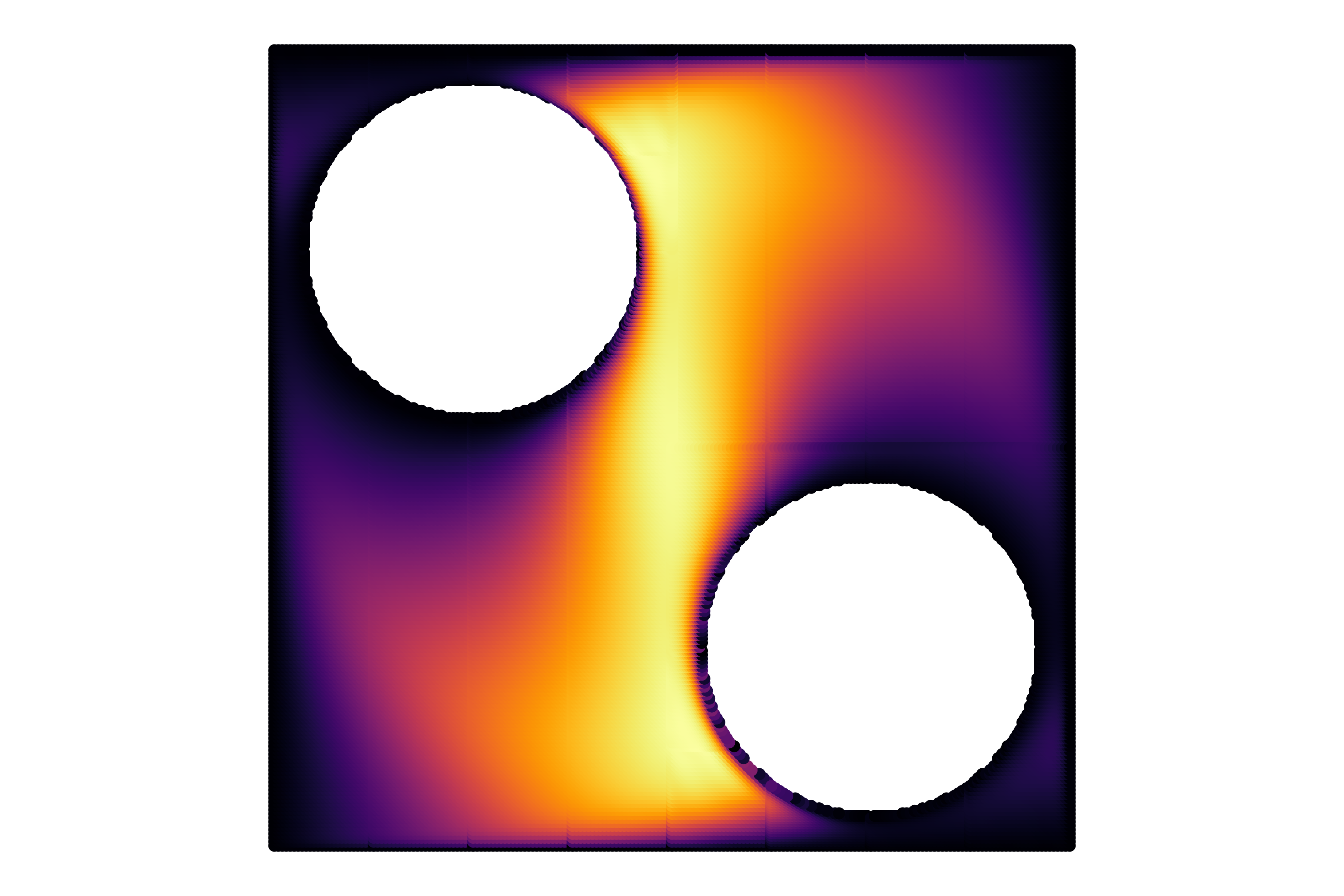}
\includegraphics[width=0.3\textwidth, trim={60em 2em 60em 2em},clip]{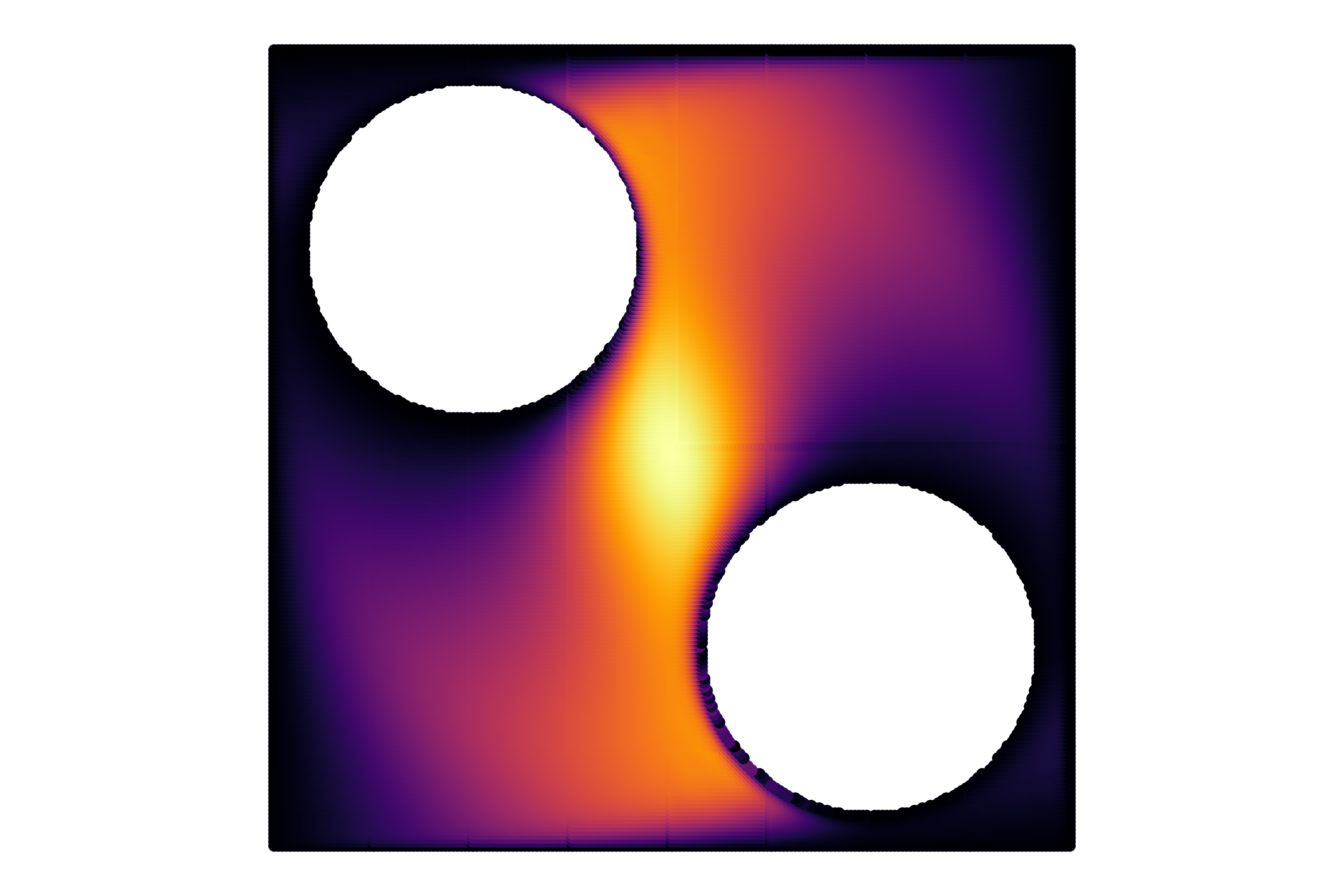}
\caption{Snapshots of the degree $P=7$ cut-cell DG solution to the advection-diffusion problem at times $t \approx 0.551$ (left), $t\approx 1.7755$ (center), and $t=3$ (right) with color limits of (0, 0.5), (0, 1.5), and (0, 2.25) respectively.}
\label{fig:cut_solution}
\end{figure}
Finally, \Cref{fig:cut_solution} shows snapshots of a degree $P=7$ solution of the advection-diffusion equation \eqref{eq:advec_diff} with zero Dirichlet boundary conditions. The advective portion is computed using a standard DG weak formulation with an upwind flux \cite{hesthaven2008nodal}, and the diffusive contribution is discretized using a Bassi-Rebay 1 viscous discretization \cite{gassner2018br1,BASSI1997267}. Zero inflow conditions are imposed on the advective discretization, with no boundary condition imposed on the outflow boundary. Zero Dirichlet boundary conditions are imposed everywhere for the viscous discretization. 

\section{Conclusion}
We have proposed a new, computationally storage- and complexity-efficient algorithm, the Givens streaming version of the \Caratheodory-Steinitz algorithm (\gscsp). We have provided mathematical stability in the total variation distance for this algorithm, and have numerically investigated the procedure and its theoretical stability on several test cases, including for pruning billion-point quadrature rules, and for generating non-standard quadrature rules in finite element simulations on non-trivial geometries. Compared to popular alternatives, the \gscsp{} algorithm is competitively stable and efficient, and requires considerably less memory.

The developments of this manuscript motivate a few directions for future work. We expect direct applicability to other PDE applications, such as kinetics where particle compression with linear moment preservation is required. Another direction could investigate whether a similar procedure could be developed to compute positive quadrature rules with upper bound constraints on the weights. Finally, while this procedure inherently addresses compression of discrete quadrature rules, an extension to quadrature rules with infinite support would significantly broaden the applicability of these procedures.

\textbf{Reproducibility of computational results.} The \gscsp{} and other related pruning algorithms were implemented in the open-source software package \texttt{CaratheodoryPruning.jl} developed by the authors, available at \url{https://github.com/fbelik/CaratheodoryPruning.jl}. A Python API for the package is available at \url{https://github.com/fbelik/CSPy}. Additionally, code for replicating the figures can be found at \url{https://github.com/fbelik/CaratheodoryFigures}.

\textbf{Acknowledgments.} FB and AN were partially supported by the Air Force Office of Scientific Research under award FA9550-23-1-0749. JC acknowledges support from the National Science Foundation under award DMS-1943186. 

\textbf{LLM Usage Declaration.} Large language models (LLMs) were used for proofreading of the manuscript. Additionally, LLMs identified minor technical errors in the appendix and supplemental materials and assisted in forming corrections. All corrections were checked and verified by the authors.

\printbibliography

\appendix

\section{Proof of \texorpdfstring{\Cref{thm:muM-perturbation}}{Stability Theorem}: Stability of \texorpdfstring{\scsp{} and \gscsp{}}{SCSP and GSCSP}}\label{sec:scsp-proof}

Throughout, we assume Assumption \ref{assum:scsp}. As discussed after Assumption \ref{assum:scsp}, \scsp{} applied to $(\mu_M,V,\Sigma)$ then behaves uniquely: every iteration has $|S| = N+1$, the cokernel vector $\bs{n}$ of $\bs{V}_{S\ast}$ is unique up to scaling, and the minimizer in \eqref{eq:sigselect-example} is unique, so exactly one node is pruned per iteration and the algorithm takes exactly $M-N$ iterations. The set of valid perturbations to $\mu_M$ is given in \Cref{def:muM-perturbations}.

At iteration $j \in [M-N]$ of \texttt{SCSP} operating on $\mu_M$, we use the following notation to identify unique objects at iteration $j$:
\begin{itemize}
  \item $\bs{n}_j$ is the cokernel vector identified in line \ref{lst:scsp-kernel-vector} of \Cref{alg:scsp}.
  \item $S_j$ is the size-$(N+1)$ set of ordered global indices in $[M]$ corresponding to the active weights at iteration $j$; $S_p \subset [M]$ denotes the size-$N$ set of indices remaining after the final iteration, i.e., the indices of $\mathrm{supp}(\nu)$.
  \item $m_j \in [N+1]$, associated with iteration $j$ of the \scsp{} algorithm, is the iteration-local index that is zeroed out.
  \item $\bs{w}_j \in \R^{N+1}$ is the $S_j$-indexed weight vector at the start of iteration $j$.
  \item $c_j \in \R$ is the constant identified by \eqref{eq:sigselect-example} such that $\bs{w}_j - c_j \bs{n}_j$ zeros out one element of $\bs{w}_j$.
\end{itemize}
There are analogous quantities arising from running \texttt{SCSP} on $\widetilde{\mu}_M$. We denote these corresponding quantities $\widetilde{\bs{n}}_j$, $\widetilde{S}_j$, $\widetilde{m}_j$, $\widetilde{\bs{w}}_j$, and $\widetilde{c}_j$, respectively.

Because the \texttt{SigSelect} function is chosen as in \eqref{eq:sigselect-example} by assumption, at iteration $j$ of the \scsp{} algorithm, the index $m_j$ and constant $c_j$ are chosen by the formulas,
\begin{align}\label{eq:analysis-m-select}
  m_j &= \argmin_{m \in S_+ \cup S_-} \left| \frac{w_m}{n_m} \right|, & c_j &= \frac{w_{m_j}}{n_{m_j}}.
\end{align}

We require some quantities defined in terms of the sequence of (unique) cokernel vectors:
\begin{align*}
  \epsilon_j &\coloneqq \min_{k \in [N+1]\backslash \{m_j\}} \left(w_{j,k} - \left|\frac{n_{j,k}}{n_{j,m_j}} \right| w_{j,m_j}\right) > 0, & N_j &\coloneqq \frac{\left\|\bs{n}_j\right\|_1}{|n_{j,m_j}|} < \infty,
\end{align*}
where $\epsilon_j > 0$ because by definition of $m_j$,
\begin{align*}
  \frac{w_{j,m_j}}{\left| n_{j,m_j} \right|} < \frac{w_{j,k}}{\left| n_{j,k} \right|} \hskip 10pt \forall \;\; k \in [N+1]\backslash\{m_j\},
\end{align*}
and $N_j < \infty$ because $n_{j,m_j}$ is non-zero. Moreover, every entry of $\bs{n}_j$ is nonzero: if $n_{j,k} = 0$ for some $k$, then the remaining $N$ entries of $\bs{n}_j$ would form a nonzero cokernel vector of an $N \times N$ submatrix of $\bs{V}$, which is invertible by the Chebyshev system assumption. In particular, $S_+ \cup S_- = [N+1]$ in \eqref{eq:analysis-m-select}, and all ratios $w_{j,k}/|n_{j,k}|$ above are finite. Note also that both $\epsilon_j$ and $N_j$ are invariant under multiplicative scaling of $\bs{n}_j$, and hence are unique numbers. The number $\epsilon_j$ is a scaled version of the optimality gap of the minimization problem \eqref{eq:analysis-m-select}, and $N_j$ measures the total mass of $\bs{n}_j$ relative to the mass on the pruned index. A final quantity we will need is a geometrically growing sequence derived from the $N_j$ numbers:
\begin{align}\label{eq:Cj}
  C_0 &= 1, & C_j &= (1 + N_j) C_{j-1} + 1, & j &> 0.
\end{align}

\begin{proof}[Proof of \Cref{thm:muM-perturbation}]
  We first define $\delta_0$. The measure $\nu = \scsp(\mu_M,V,\Sigma)$ is unique, having support points and weights,
  \begin{align}\label{eq:nu-def}
    \nu &= \sum_{\ell \in [N]} u_\ell \delta_{z_\ell}, & \{z_{\ell}\}_{\ell \in [N]} &\subset \mathrm{supp}(\mu_M), & \underbar{u} &\coloneqq \min_{\ell \in [N]} u_\ell > 0,
  \end{align}
  where $\underbar{u} > 0$ because Assumption \ref{assum:scsp} guarantees that each step of the \scsp{} algorithm zeros out exactly one weight. With $\bs{V} \in \R^{M \times N}$ the Vandermonde-like matrix defined in \eqref{eq:moment-conditions}, and $S_p$ as above, define
  \begin{align}\label{eq:UD-def}
    \bs{U} &\coloneqq \bs{V}_{S_p\ast} \bs{D}^{-1}, & \bs{D} &= \mathrm{diag}\left(\|v_1\|_{L_\mu^1(X)},\, \ldots,\, \|v_N\|_{L_\mu^1(X)} \right),
  \end{align}
  and note that $\bs{V}_{S_p\ast}$ is invertible by the Chebyshev system assumption. We also abbreviate
  \begin{align}\label{eq:kappa-def}
    \kappa \coloneqq \frac{\sqrt{N}\, \|\bs{U}^{-1}\|_2}{\tau}.
  \end{align}
  Then we define $\delta_0$ as,
  \begin{align*}
    \delta_0 &\coloneqq \min_{d \in [5]}\left\{\delta_d\right\}, &
    \delta_1 &= \frac{1}{3}, &
    \delta_2 &= \frac{1}{3|\mu_M|} \min_{j \in [M-N]} \frac{\epsilon_j}{C_j}, \\
    \delta_3 &= \frac{\underbar{u}}{12 \kappa |\mu_M|}, &
    \delta_4 &= \frac{\underbar{u}}{6 |\mu_M| C_{M-N}}, &
    \delta_5 &= \frac{|\nu|}{6 |\mu_M| (C_{M-N} + N\kappa)}.
  \end{align*}
  Now consider any $(\widetilde{\mu}_M, \widetilde{\Sigma}) \in P_\tau(\mu_M,\Sigma)$ satisfying $\delta = d_{\mathrm{TV}}(\mu_M, \widetilde{\mu}_M) < \delta_0$. If $\delta = 0$, then $(\widetilde{\mu}_M,\widetilde{\Sigma}) = (\mu_M,\Sigma)$ and there is nothing to prove, so we assume $\delta > 0$. Let $\widetilde{M} = |\mathrm{supp}(\widetilde{\mu}_M)| \geq M$ denote the support size of $\widetilde{\mu}_M$, and let $\bs{w}, \widetilde{\bs{w}} \in \R^{\widetilde{M}}$ denote the weights on these support points, ordered according to $\widetilde{\Sigma}$. Since $\widetilde{\Sigma}_{[M]} = \Sigma$, entries $1,\ldots,M$ correspond to $\mathrm{supp}(\mu_M)$ in the order $\Sigma$, and entries $M+1,\ldots,\widetilde{M}$ correspond to the nodes added by the perturbation.
  The $\widetilde{M}$-vector $\bs{w}$ is formed by padding the original weights $\bs{w}_{[M]}$ for $\mu_M$ with $\widetilde{M} - M$ zeros. I.e.,
  \begin{align*}
    \widetilde{\bs{w}} &= \big( \widetilde{w}_1, \ldots, \widetilde{w}_M, \widetilde{w}_{M+1}, \ldots, \widetilde{w}_{\widetilde{M}} \big), & \bs{w} &= \big( w_1, \ldots, w_M, \underbrace{0, \ldots, 0}_{\widetilde{M}-M \textrm{ entries}} \big).
  \end{align*}
  Note that this zero padding of $\bs{w}$ does not affect the result of $\nu = \scsp(\mu_M,V,\Sigma)$ since after $M-N$ iterations the procedure would simply prune the zero-padded weights because the Chebyshev system assumption ensures that $n_{j,N+1} \neq 0$. 
  With this setup, $d_{\mathrm{TV}}(\mu_M, \widetilde{\mu}_M) = \delta$ means $\|\bs w - \widetilde{\bs w}\|_1 = \delta(\|\bs w\|_1 + \|\widetilde{\bs w}\|_1)$. By the reverse triangle inequality and $\delta < \delta_1 = 1/3$, $\|\widetilde{\bs w}\|_1 - \|\bs w\|_1 < \frac{1}{3}(\|\bs w\|_1 + \|\widetilde{\bs w}\|_1)$, i.e., $\|\widetilde{\bs w}\|_1 < 2\|\bs w\|_1$. Hence, for any $T \subset [\widetilde{M}]$,
  \begin{align}\label{eq:wT-bound}
    \left\| \bs{w}_T - \widetilde{\bs{w}}_T \right\|_1 \leq \|\bs w - \widetilde{\bs w}\|_1 = \delta\left(\|\bs w\|_1 + \|\widetilde{\bs w}\|_1\right) < 3 \delta \|\bs{w}\|_1 = 3 \delta |\mu_M|.
  \end{align}
  In particular,
  \begin{align}\label{eq:wtilde-1}
    \sum_{i=M+1}^{\widetilde{M}} \widetilde{w}_i &< 3 \delta |\mu_M|, &
    |w_i - \widetilde{w}_i| &< 3\delta|\mu_M|, & i &\in [M].
  \end{align}

  We now analyze $\widetilde{\nu} = \scsp(\widetilde{\mu}_M,V,\widetilde{\Sigma})$, which we break up into two parts. The first part considers the first $M-N$ iterations, which operate on $\bs{w}_{[M]}$ and $\widetilde{\bs{w}}_{[M]}$ that involve nodes in the shared set $\mathrm{supp}(\mu_M)$. The second part of the analysis considers nodes in the set $\mathrm{supp}(\widetilde{\mu}_M) \backslash \mathrm{supp}(\mu_M)$ that are supported only in $\widetilde{\mu}_M$.

  For the first part of the analysis, fix $j \in [M-N]$ and make the inductive hypothesis that at the start of iteration $j$,
  \begin{align}\label{eq:inductive-hypothesis}
    \left\|\bs{w}_j - \widetilde{\bs{w}}_j \right\|_1 &< 3\delta\left|\mu_M\right| C_{j-1} , & S_j &= \widetilde{S}_j.
  \end{align}
  For $j = 1$ this holds with $C_0 = 1$: both algorithms start from the first $N+1$ indices, and the bound is \eqref{eq:wT-bound} with $T = [N+1]$. Since $\delta < \delta_2$ and, by \eqref{eq:Cj}, $C_{j-1}/C_j = (1 - C_j^{-1})/(1+N_j) < (1+N_j)^{-1}$,
  \begin{align}\label{eq:w-eps-strict-inequality}
    \left\|\bs{w}_j - \widetilde{\bs{w}}_j \right\|_1 < 3\left|\mu_M\right| \delta_2 C_{j-1} \leq \epsilon_j \frac{C_{j-1}}{C_j} < \frac{\epsilon_j}{1 + N_j}.
  \end{align}

  Since $\widetilde{S}_j = S_j$, the perturbed algorithm computes a cokernel vector of the same matrix $\bs{V}_{S_j\ast}$, so $\widetilde{\bs{n}}_j$ is a nonzero multiple of $\bs{n}_j$. The pruned vector $\widetilde{\bs{w}}_j - \widetilde{c}_j \widetilde{\bs{n}}_j$ is unchanged by rescaling $\widetilde{\bs{n}}_j$, so we may take $\widetilde{\bs{n}}_j = \bs{n}_j$. 
  Since $|n_{j,k}| \leq \|\bs{n}_j\|_1$, we have $|n_{j,k}/n_{j,m_j}| \leq N_j$. Then for any $k \neq m_j$,
  \begin{align*}
    w_{j,k} - \left|\frac{n_{j,k}}{n_{j,m_j}}\right|\, w_{j,m_j} &\geq \epsilon_j \\ 
    &\stackrel{\eqref{eq:w-eps-strict-inequality}}{>} (1+N_j)\|\widetilde{\bs{w}}_j - \bs{w}_j\|_1 \\
    &\geq N_j \left|\widetilde{w}_{j,m_j} - w_{j,m_j}\right| + \left|\widetilde{w}_{j,k} - w_{j,k}\right| \\
    &\geq \left|\frac{n_{j,k}}{n_{j,m_j}}\right| \left(\widetilde{w}_{j,m_j} - w_{j,m_j}\right) + \left(w_{j,k} - \widetilde{w}_{j,k}\right).
  \end{align*}
  Comparing the first and last expressions yields
  \begin{align*}
    \frac{\widetilde{w}_{j,m_j}}{\left| n_{j,m_j}\right|} &< \frac{\widetilde{w}_{j,k}}{\left| n_{j,k}\right|}, & k &\in [N+1]\backslash\{m_j\}.
  \end{align*}
  Hence, the perturbed version of the minimization problem \eqref{eq:analysis-m-select} identifies the same index as the unperturbed problem, $\widetilde{m}_j = m_j$, i.e., the same node is chosen for removal in both algorithms. Since both algorithms then append the node at the same global index $N+j+1$, this also gives $S_{j+1} = \widetilde{S}_{j+1}$. The two step sizes are $c_j = w_{j,m_j}/n_{j,m_j}$ and $\widetilde{c}_j = \widetilde{w}_{j,m_j}/n_{j,m_j}$, so $|c_j - \widetilde{c}_j| = |\widetilde{w}_{j,m_j} - w_{j,m_j}|/|n_{j,m_j}|$, and the difference between the corresponding iteration-$j$ pruned weight vectors is,
  \begin{align*}
    \left\| \left( \bs{w}_j - c_j \bs{n}_j \right) - \left( \widetilde{\bs{w}}_j - \widetilde{c}_j \bs{n}_j \right)\right\|_1 &\leq \left\| \widetilde{\bs{w}}_j - \bs{w}_j \right\|_1 + N_j |\widetilde{w}_{j,m_j} - w_{j,m_j}| \\
    &\leq \left\| \widetilde{\bs{w}}_j - \bs{w}_j \right\|_1 (1 + N_j) \\
    &\stackrel{\eqref{eq:inductive-hypothesis}}{<} 3 \delta |\mu_M| C_{j-1} (1 + N_j).
  \end{align*}
  The algorithm then forms $\bs{w}_{j+1}$ and $\widetilde{\bs{w}}_{j+1}$ by deleting entry $m_j$, which is zero in both pruned vectors and so does not change the $\ell^1$ norm of their difference, and appending the original weights at global index $N+j+1$. Hence, by \eqref{eq:wtilde-1},
  \begin{align*}
    \left\| \bs{w}_{j+1} - \widetilde{\bs{w}}_{j+1} \right\|_1 &= \left| w_{N+j+1} - \widetilde{w}_{N+j+1} \right| + \left\| \left( \bs{w}_j - c_j \bs{n}_j \right) - \left( \widetilde{\bs{w}}_j - \widetilde{c}_j \bs{n}_j \right)\right\|_1 \\
    &< 3\delta|\mu_M| \left(1 + C_{j-1} (1 + N_j)\right) = 3 \delta |\mu_M| C_j,
  \end{align*}
  which completes the proof of \eqref{eq:inductive-hypothesis} for iteration $j+1$. By finite induction, we conclude that after $M-N$ iterations have completed, we have pruned weight vectors $\bs{w}_{M-N+1}$ and $\widetilde{\bs{w}}_{M-N+1}$, which satisfy,
  \begin{align}\label{eq:C-bound-1}
    \left\| \bs{w}_{M-N+1} - \widetilde{\bs{w}}_{M-N+1} \right\|_1 
    < 3 \delta |\mu_M| C_{M-N}.
  \end{align}
  Since $w_{M-N+1,k} \geq \underbar{u}$ for $k \in [N]$ by \eqref{eq:nu-def}, inserting $\delta < \delta_4$ above yields
  \begin{equation}\label{eq:w-first-lb}
    \left\| \bs{w}_{M-N+1} - \widetilde{\bs{w}}_{M-N+1} \right\|_1 < \frac{1}{2}\underbar{u} \implies \widetilde{w}_{M-N+1,k} > \frac{1}{2}\underbar{u}\;\;\forall\;\;k\in[N].
  \end{equation}

  For the second part of the analysis, we consider iterations $j \in [M-N+1, \widetilde{M}-N]$ (there are none if $\widetilde{M} = M$). Only $\scsp(\widetilde{\mu}_M,V,\widetilde{\Sigma})$ is still running, so we write $S_j$ and $\bs{n}_j$ for $\widetilde{S}_j$ and $\widetilde{\bs{n}}_j$, and order $S_j$ so that its first $N$ entries are $S_p$. Let $z_1,\ldots,z_N$ be the nodes indexed by $S_p$, as in \eqref{eq:nu-def}. Our inductive hypothesis is that the only node outside $S_p$ is the one most recently appended, and that the retained weights stay bounded away from zero:
  \begin{align}\label{eq:inductive}
    S_j &= S_p \cup \{j+N\}, & \min_{q \in [N]} \widetilde{w}_{j,q} &> \frac{\underbar{u}}{4}.
  \end{align}
  For $j = M-N+1$, the first relation holds because both algorithms prune the same indices in the first part, and the second is \eqref{eq:w-first-lb}.

  Let $x \coloneqq x_{j+N} \in X_\tau$ be the node at global index $j+N$ (local index $N+1$). Since $\bs{V}_{S_p\ast}$ is invertible, $\bs{V}_{S_j\ast}$ has a one-dimensional cokernel, and its cokernel vector satisfies $n_{j,N+1} \neq 0$ and
  \begin{align}\nonumber
    \bs{0} = (\bs{V}_{S_j\ast})^T \bs n_j &= (\bs{V}_{S_p\ast})^T \bs{n}_{j,[N]} + n_{j,N+1}\, \bs{v}(x) \\ \label{eq:n-rearranged}
    \implies \quad \frac{\bs{n}_{j,[N]}}{n_{j,N+1}} &= -\left(\bs{V}_{S_p\ast}\right)^{-T} \bs{v}(x) = -\bs{U}^{-T}\bs{D}^{-1}\bs{v}(x).
  \end{align}
  Each entry of $\bs{D}^{-1}\bs{v}(x)$ is $v_q(x)/\|v_q\|_{L_\mu^1(X)}$, of magnitude at most $1/\tau$ since $x \in X_\tau$. Using $\|\bs{A}\|_\infty \leq \sqrt{N}\|\bs{A}\|_2$ for $\bs{A} \in \R^{N\times N}$ and $\|\bs{U}^{-T}\|_2 = \|\bs{U}^{-1}\|_2$,
  \begin{align}\label{eq:n-frac-simplified}
    \frac{\|\bs{n}_{j,[N]}\|_{\infty}}{|n_{j,N+1}|} \leq \| \bs{U}^{-T} \|_\infty \| \bs{D}^{-1} \bs{v}(x)\|_\infty \leq \frac{\sqrt{N} \|\bs{U}^{-1}\|_2}{\tau} = \kappa.
  \end{align}

  The appended weight has not been modified, so $\widetilde{w}_{j,N+1} = \widetilde{w}_{j+N} < 3\delta|\mu_M|$ by \eqref{eq:wtilde-1}. Thus, for any $q \in [N]$,
  \begin{align*}
    \frac{|n_{j,q}|}{|n_{j,N+1}|} \widetilde{w}_{j,N+1} \stackrel{\eqref{eq:n-frac-simplified}}{\leq} \kappa\, \widetilde{w}_{j,N+1} < 3 \kappa \delta |\mu_M| \stackrel{(\delta<\delta_3)}{<} \frac{\underbar{u}}{4} \stackrel{\eqref{eq:inductive}}{<} \widetilde{w}_{j,q}.
  \end{align*}
  Dividing by $|n_{j,q}|$ when $n_{j,q} \neq 0$ (otherwise $q$ is not a candidate in \eqref{eq:analysis-m-select}) shows that local index $N+1$ is the unique minimizer in \eqref{eq:analysis-m-select}. Hence, $x$ is removed with $\widetilde{c}_j = \widetilde{w}_{j,N+1}/n_{j,N+1}$, the next node is appended so that $S_{j+1} = S_p \cup \{j+N+1\}$, and $\widetilde{\bs{w}}_{j+1,[N]} = \widetilde{\bs{w}}_{j,[N]} - \widetilde{c}_j \bs{n}_{j,[N]}$. Unrolling this update back to iteration $M-N+1$ gives
  \begin{equation}\label{eq:unravel-w}
    \widetilde{\bs{w}}_{j+1,[N]} = \widetilde{\bs{w}}_{M-N+1,[N]} - \sum_{\ell=M-N+1}^{j} \frac{\widetilde{w}_{\ell,N+1}}{n_{\ell,N+1}} \bs{n}_{\ell,[N]}.
  \end{equation}

  For $q \in [N]$, using \eqref{eq:wtilde-1}, \eqref{eq:n-frac-simplified}, and \eqref{eq:unravel-w},
  \begin{align}\nonumber
    \widetilde{w}_{j+1,q} \geq\; &\widetilde{w}_{M-N+1,q} - \kappa \sum_{\ell=M-N+1}^{j} \widetilde{w}_{\ell,N+1} > \widetilde{w}_{M-N+1,q} - 3\kappa\delta|\mu_M| \\ \label{eq:w-lb-derivation}
    \stackrel{(\delta < \delta_3)}{>} &\widetilde{w}_{M-N+1,q} - \frac{\underbar{u}}{4} \stackrel{\eqref{eq:w-first-lb}}{>} \frac{\underbar{u}}{4}.
  \end{align}
  This establishes \eqref{eq:inductive} for iteration $j+1$.

  After the final iteration $j = \widetilde{M}-N$, the active set is $S_p$, so $\widetilde{\nu}$ is supported on $z_1,\ldots,z_N$, the same nodes as $\nu$. Applying \eqref{eq:unravel-w} with $j = \widetilde{M}-N$, the triangle inequality, $\|\bs{y}\|_1 \leq N\|\bs{y}\|_\infty$ for $\bs{y} \in \R^N$, \eqref{eq:wtilde-1}, and \eqref{eq:n-frac-simplified},
  \begin{align}\nonumber
    \left\| \widetilde{\bs{w}}_{\widetilde{M}-N+1,[N]} - \widetilde{\bs{w}}_{M-N+1,[N]} \right\|_1 &\leq N \sum_{\ell=M-N+1}^{\widetilde{M}-N} \widetilde{w}_{\ell,N+1} \frac{\|\bs{n}_{\ell,[N]}\|_\infty}{|n_{\ell,N+1}|} \\ \label{eq:C-bound-2}
    &\leq N\kappa \sum_{i=M+1}^{\widetilde{M}} \widetilde{w}_i < 3 \delta |\mu_M| N\kappa.
  \end{align}
  Combining \eqref{eq:C-bound-1} and \eqref{eq:C-bound-2} with the triangle inequality, and then using $\delta < \delta_5$,
  \begin{align*}
    \left\|\bs{w}_{M-N+1,[N]} - \widetilde{\bs{w}}_{\widetilde{M}-N+1,[N]} \right\|_1 < 3\delta |\mu_M| ( C_{M-N} + N\kappa ) < \frac{|\nu|}{2}.
  \end{align*}
  By the reverse triangle inequality, $|\widetilde{\nu}| > |\nu| - |\nu|/2 = |\nu|/2$, so $|\nu| + |\widetilde{\nu}| > \frac{3}{2}|\nu|$. Since $\nu$ and $\widetilde{\nu}$ share their support,
  \begin{align*}
    d_{\mathrm{TV}}(\nu,\widetilde{\nu}) &= \frac{\left\|\bs{w}_{M-N+1,[N]} - \widetilde{\bs{w}}_{\widetilde{M}-N+1,[N]} \right\|_1}{|\nu| + |\widetilde{\nu}|} < C \delta, &
    C &= \frac{2|\mu_M|}{|\nu|} (C_{M-N} + N\kappa).
  \end{align*}
\end{proof}

\section{\Caratheodory-Steinitz pruning visualization}\label{sec:supp-cs-pruning}
This section provides the visual depiction of the two pruning choices cited in the main text.

\begin{figure}[htbp]
    \centering
    \includegraphics[width=0.99\textwidth]{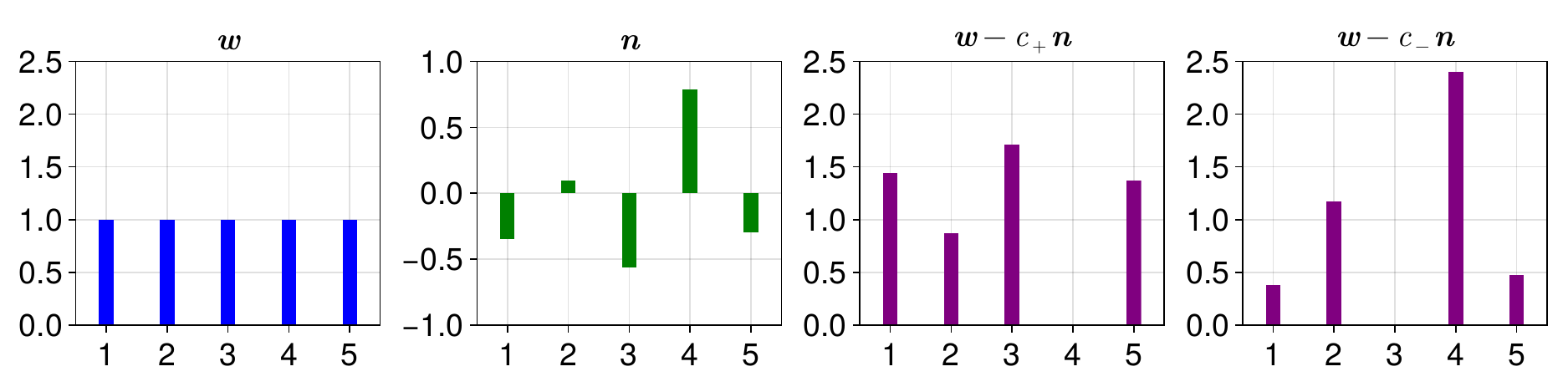}
    \caption{Visual depiction of the two possible pruning choices for given weights and kernel vector. From left to right: visualizing $\bs w$, $\bs n$, $\bs w - c_{+}\bs n$, and $\bs w - c_{-}\bs n$.}
    \label{fig:cs-pruning}
\end{figure}

\section{Tensorized quadrature attaining \texorpdfstring{$Q < N$}{Q < N}}
\label{sec:gq-highd}
This section provides explicit examples of Tchakaloff quadrature rules satisfying \eqref{eq:tchakaloff} with $Q < N$. Hence, while in this paper we consider identifying such rules with $Q = N$, the example in this section demonstrates that such rules need not be minimal quadrature rules. Fix $k \in \N_0$, $d \in \N$, and let $\mu$ be a product measure on $X=\R^d$, and let $V$ be the subspace of at-most degree-$k$ $d$-variate polynomials. With $\mu_j$ the coordinate-$j$ marginal measure of $\mu$, assume $\mu_j$ has finite moments up to univariate degree $k+1$, and let $\{\bs{x}^{(j)}_q, w^{(j)}_q\}_{q \in [p]}$ be the $p$-point univariate $\mu_j$-Gaussian quadrature rule, where $p \coloneqq \left\lceil \frac{k+1}{2} \right\rceil$. This choice of $p$ ensures exact $\mu_j$-integration of degree $2p - 1 \geq k$ polynomials. By tensorizing these $d$ different $p$-point rules, we have a $Q = p^d$-point quadrature rule that exactly $\mu$-integrates all $d$-variate polynomials of degree at most $k$. We can now compare $Q$ and $N = \dim(V)$:
\begin{align*}
N &= \left(\begin{array}{c} k+d \\ d \end{array}\right) = \frac{(k+1)^{(d)}}{d!} \sim \frac{k^d}{d!}, &
  \frac{Q}{N} &=\frac{\left\lceil \frac{k+1}{2} \right\rceil^d}{(k+1)^{(d)}} d! \stackrel{k \gg 1}{\sim} \frac{d!}{2^d}, 
\end{align*}
where $(k+1)^{(d)}$ denotes the rising factorial/Pochhammer function, $(k+1)^{(d)} = \prod_{j=1}^{d} (k+j)$. Hence, for large $k, d$, the above ratio is greater than unity, implying $Q > N$, so that this is not a Tchakaloff-attaining quadrature. However, direct computation of $Q/N$ when $k = 2$ shows that this construction achieves $Q < N$ when $d < 4$. 

\section{Givens rotation-based downdates and updates}\label{sec:gud}
This section provides more detailed pseudocode, \Cref{alg:gud}, that accomplishes the procedure described in \Cref{ssec:gscsp}: An $\mathcal{O}(N^2)$ procedure that uses Givens rotations to update the full QR decomposition of an $(N+k) \times N$ matrix by replacing one row of the original matrix with a new row. The \gscsp{} algorithm is the \scsp{} procedure in \Cref{alg:scsp} augmented by using \Cref{alg:gud} as a subroutine to accomplish line \ref{lst:scsp-augment} of \Cref{alg:scsp}.

\begin{algorithm}[htbp]
    \caption{Givens Row UpDowndate and Kernel Vector}\label{alg:gud}
    \begin{flushleft}
        \textbf{Input: } $\bs{V} \in \mathbb{R}^{M\times N}$, $T \subset [M]$ with $|T| = N+k$, $\bs{Q}\in\mathbb{R}^{N+k\times N+k}$, $\bs{R}\in\mathbb{R}^{N+k\times N}$, $i_\text{rem}\in T$, $i_\text{new} \in [M] \setminus T$\\
        \textbf{Output: } $\bs{Q} \in\mathbb{R}^{N+k\times N+k}$, $\bs{R}\in\mathbb{R}^{N+k\times N}$, $\bs{n}\in\R^{N+k}$, $T$\\
    \end{flushleft}
    \begin{algorithmic}[1]
        \State $j_\text{rem} \gets \text{indexof}(T, i_\text{rem})$
        \State $T \gets (T \setminus \{i_\text{rem}\}) \cup \{i_\text{new}\}$
        \For{$i=N+k$ down to $j_\text{rem}+1$} \Comment{Begin Givens downdate}
            \State Form Givens rotation $\bs{G}$ for indices $i$ and $i-1$ such that $(\bs{Q}\bs{G}^T)_{j_\text{rem},i} = 0$ \Comment{$O(1)$}
            \State $\bs{Q} \gets \bs{Q}\bs{G}^T$ \Comment{$O(N+k)$, repeated $N+k-j_\text{rem}$ times}
            \State $\bs{R} \gets \bs{G}\bs{R}$ \Comment{$O(N)$, repeated $N+k-j_\text{rem}$ times}
        \EndFor
        \For{$i=j_\text{rem}-1$ down to $1$}
            \State Form Givens rotation $\bs{G}$ for indices $i$ and $j_\text{rem}$ such that $(\bs{Q}\bs{G}^T)_{j_\text{rem},i} = 0$ \Comment{$O(1)$}
            \State $\bs{Q} \gets \bs{Q}\bs{G}^T$ \Comment{$O(N+k)$, repeated $j_\text{rem}-1$ times}
            \State $\bs{R} \gets \bs{G}\bs{R}$ \Comment{$O(N)$, repeated $j_\text{rem}-1$ times}
        \EndFor \Comment{End Givens downdate}
        \State $\bs{R}_{\{j_\text{rem}\}*} \gets \bs{V}_{\{i_\text{new}\}*}$ \Comment{Begin Givens update}
        \State $\bs{Q}_{j_\text{rem},j_\text{rem}} \gets +1.0$
        \For{$i=1$ to $\min(N,j_\text{rem}-1)$}
            \State Form Givens rotation $\bs{G}$ for indices $i$ and $j_\text{rem}$ such that $(\bs{G}\bs{R})_{j_\text{rem},i} = 0$ \Comment{$O(1)$}
            \State $\bs{Q} \gets \bs{Q}\bs{G}^T$ \Comment{$O(N+k)$, repeated $\min(N,j_\text{rem}-1)$ times}
            \State $\bs{R} \gets \bs{G}\bs{R}$ \Comment{$O(N)$, repeated $\min(N,j_\text{rem}-1)$ times}
        \EndFor
        \For{$i=j_\text{rem}$ to $N$}
            \State Form Givens rotation $\bs{G}$ for indices $i+1$ and $i$ such that $(\bs{G}\bs{R})_{i+1,i} = 0$ \Comment{$O(1)$}
            \State $\bs{Q} \gets \bs{Q}\bs{G}^T$ \Comment{$O(N+k)$, repeated $\max(0,N-j_\text{rem}+1)$ times}
            \State $\bs{R} \gets \bs{G}\bs{R}$ \Comment{$O(N)$, repeated $\max(0,N-j_\text{rem}+1)$ times}
        \EndFor \Comment{End Givens update}
        \State $\bs{n} \gets \bs{Q}_{*\{N+1\}}$
    \end{algorithmic}
\end{algorithm}

\section{Stability of NNLS}
\label{app:nnls_stability}

A simplified Lawson-Hanson algorithm, without details of efficient updates and downdates by Householder reflections, is provided in \Cref{alg:nnls} \cite{lawson_1995,bro1997fast}. The algorithm iteratively constructs a more accurate and less-sparse solution, $\bs{w}$, by including indices with the largest dual and by solving least squares problems. The dual at a given iteration is defined proportional to the gradient of the squared $\ell^2$ moment error:
\[\bs{d} = \bs{V}(\bs{\eta} - \bs{V}^T\bs{w}) = \frac{-1}{2}\nabla_{\bs{w}} \|\bs{V}^T\bs{w} - \bs{\eta}\|_2^2.\]
For simplicity of analysis, \Cref{alg:nnls} is a simplified variant of the Lawson-Hanson method: in the inner loop, the index attaining the minimum ratio is removed, but the iterate $\bs{w}$ is not moved toward $\bs{s}$ before recomputing the least squares solution. If this variant terminates, its output satisfies the Karush-Kuhn-Tucker conditions and hence solves the NNLS problem; the stability analysis below applies to this variant. We note at the end of the proof that this can be extended to the full NNLS algorithm, albeit with more complicated constants.
In the algorithm, an index set $P \subset [M]$ is gradually built, and the final solution satisfies nonnegativity, zero error gradient for indices with positive weight, and positive error gradient for active indices (zero weights). The inner loop (starting on line \ref{line:innerloop} of \cref{alg:nnls}) is intended to occur infrequently and plays the role of removing indices whose weights are made negative by adding new indices.
\begin{algorithm}[H]
    \caption{\nnls{}: Simplified Lawson-Hanson NNLS Algorithm}\label{alg:nnls}
    \begin{flushleft}
        \textbf{Input: } $\bs{V} \in \mathbb{R}^{M\times N}$, $\bs{\eta} = \bs{V}^T \bs{w} \in \R^N$\\
        \textbf{Output: } $P \subset [M]$ and $\bs{w} \in \R^{M}_+$ which is a $P$-sparse vector that solves $\min_{\bs{w}\geq0}\|\bs{V}^T\bs{w}-\bs{\eta}\|_2$\\
    \end{flushleft}
    \begin{algorithmic}[1]
        \State $P \gets \emptyset$, $\bs{w}\gets\bs{0}$, $\bs{s}\gets\bs{0}$
        \State $\bs{d} \gets \bs{V}(\bs{\eta} - \bs{V}^T\bs{w}) = \bs{V}\bs{\eta}$
        \While{$\max(\bs{d}) > 0$} \label{line:nnlsouterloop}
            \State $m \gets \mathrm{maxindex}(\bs{d})$, $P\gets P\cup\{m\}$ \label{line:maxindex}
            \State $\bs{s}_P \gets ({\bs{V}_{P*}}^{T})^\dagger \bs{\eta}$, $\bs{s}_{[M]\backslash P} \gets \bs{0}$
            \State $Q \gets \{i\in P : s_i \leq 0\}$ \label{line:Qdef}
            \While{$|Q| > 0$} \label{line:innerloop}
                \State $i_\mathrm{rem} \gets \text{argmin}_{i\in Q} \frac{-w_i}{s_i - w_i}$ \label{line:remindex}
                \State $P \gets P \backslash \{i_\mathrm{rem}\}$, $w_{i_\mathrm{rem}} \gets 0$ \label{line:wzero}
                \State $\bs{s}_P \gets ({\bs{V}_{P*}}^{T})^\dagger \bs{\eta}$, $\bs{s}_{[M]\backslash P} \gets \bs{0}$
                \State $Q \gets \{i\in P : s_i \leq 0\}$ 
            \EndWhile
            \State $\bs{w}_P \gets \bs{s}_P$ \label{line:wupdate}
            \State $\bs{d} \gets \bs{V}(\bs{\eta} - \bs{V}^T\bs{w})$
        \EndWhile
    \end{algorithmic}
\end{algorithm}

Fixing $\mu_M \in P_+$, define the set of admissible \nnls{} perturbations as
\begin{align*}
  P_{\mathrm{NNLS}}(\mu_M) &\coloneqq \left\{ \mu \in P_+ \;\;\big|\;\; \mathrm{supp} (\mu) = \mathrm{supp} (\mu_M) \right\}.
\end{align*}
Similar to the proof of stability for the \scsp{} algorithm, we require the following assumptions on the measure $\mu_M$ and the hyperparameters of the \nnls{} algorithm.
\begin{assumption}\label{assum:nnls}
  Suppose the \nnls{} algorithm is run on $(V, \mu_M)$. We assume that:
  \begin{itemize}
    \item With $V$ fixed, then running the \nnls{} algorithm for any $\mu_M$ uses the same basis $v_n(\cdot)$ as in \eqref{eq:v-def}.
    \item The optimizations have unique solutions for all iterations (specifically lines \ref{line:maxindex} and \ref{line:remindex} of \Cref{alg:nnls}).
    \item For all iterations, the intermediate least squares solutions have full density, i.e., for a realized index set $P$, the least squares solution on those indices has exactly $|P|$ nonzero elements.
    \item The algorithm terminates, and its output has exactly $N$ nonzero elements.
  \end{itemize}
\end{assumption}
Under \Cref{assum:nnls}, the \nnls{} algorithm is robust to $d_{\mathrm{TV}}$ perturbations in $P_{\mathrm{NNLS}}$.
\begin{theorem}[\nnls{} stability]\label{thm:muM-perturbation-nnls}
  Fix $V \subset L^1(X)$, and let $\mu_M \in P_+$ be a given measure, and let $\nu = \nnls(\mu_M,V)$ be the output of \Cref{alg:nnls} satisfying \Cref{assum:nnls}. Then the \nnls{} algorithm is locally Lipschitz (in particular continuous) with respect to the total variation distance in a $d_{\mathrm{TV}}$-neighborhood of $P_{\mathrm{NNLS}}(\mu_M)$ around $\mu_M$. I.e., there are positive constants $\delta_0 = \delta_0(\mu_M)$ and $C = C(\mu_M)$ such that for any $\widetilde{\mu}_M \in P_{\mathrm{NNLS}}(\mu_M)$ satisfying $d_{\mathrm{TV}}(\mu_M, \widetilde{\mu}_M) < \delta_0$, then
\begin{align*}
  d_{\mathrm{TV}}(\nu, \widetilde{\nu}) \leq C\, d_{\mathrm{TV}}(\mu_M, \widetilde{\mu}_M),
\end{align*}
  where $\widetilde{\nu} = \nnls(\widetilde{\mu}_M,V)$.
\end{theorem}

\begin{proof}[Proof of \Cref{thm:muM-perturbation-nnls}]
  Having provided the rigorous details for the proof of \Cref{thm:muM-perturbation}, we provide here mainly a sketch, omitting many steps, technical notations, and computations that are similar in spirit for this proof. For example, we will only seek to show that $\|\bs{u} - \widetilde{\bs{u}}\|_2$ is small, with $\bs{u} \in \R^N$ the weights for $\nu$ and $\widetilde{\bs{u}} \in \R^N$ the weights for $\widetilde{\nu}$, omitting the computations that connect these quantities to the total variation distances. We will also provide the argument for a single outer loop iteration of the algorithm instead of providing the meticulous argument that holds for every iteration.

  We let $\bs{V}\in\R^{M\times N}$ and $\bs{w}\in\R^M$  be the Vandermonde matrix and weight vector associated with $(V,\mu_M)$. Recall we assume that $\bs{u}$ is the output of \nnls($\bs{V},\bs{V}^T\bs{w})$, and that $\bs{u}\in\R^N$. Let $S_1\subset \mathcal{P}([M])$ (the power set of $[M]$) be the collection of index sets built by the unperturbed \nnls{} algorithm at the start of the outer loop (the index sets $P$ that are encountered on line \ref{line:nnlsouterloop} of \Cref{alg:nnls}). Similarly, let $S_2\subset \mathcal{P}([M])$ be the set of index sets $P$ observed by the unperturbed \nnls{} algorithm on line \ref{line:innerloop} each time its condition is evaluated (including when the inner loop is not entered). Finally, let $P_0$ be the final set of indices, satisfying $|P_0|=N$ by \Cref{assum:nnls}.

  For any $P\in S_1\cup S_2$, we define $\bs{s}(P)$ to be the weight vector corresponding to that index set, given by $\bs{s}(P)_P = (\bs{V}_{P*}^T)^\dagger\bs{\eta}$ and $\bs{s}(P)_{[M]\backslash P} = \bs{0}$. We similarly define $\bs{d}(P)=\bs{V}(\bs{\eta}-\bs{V}_{P*}^T\bs{s}(P)_P)$ for any $P\in S_1$ to be the corresponding dual vector at the iteration corresponding to that value of $P$. 
  For $P \in S_2$, let $P' = P'(P) \in S_1$ be the index set at the start of the outer iteration in which $P$ is observed (i.e., before the index $m$ is added on line \ref{line:maxindex}), and let $\bs{w}(P) \in \R^M$ be the iterate $\bs{w}$ when the condition on line \ref{line:innerloop} is evaluated with index set $P$. Since $\bs{w}$ is only modified on lines \ref{line:wzero} and \ref{line:wupdate}, $w(P)_i = s(P')_i$ for $i \in P \cap P'$, and $w(P)_i = 0$ otherwise. When $P' = \emptyset$ (so that $\bs{w}(P) = \bs{0}$), we interpret $\|(\bs{V}_{P'*}^T)^\dagger\bs{V}^T\|_\infty$ as $0$.
  We also define $Q = Q(P)$ to be the set defined on line \ref{line:Qdef}. We use tildes to denote weight vectors and dual vectors in the perturbed problem.

    Our main effort will seek to ensure that the discrete optimization problems on lines \ref{line:maxindex} and \ref{line:remindex} have the same solutions for the unperturbed and the perturbed problems. To this end, we define optimality gaps $\delta_1(P)$ and $\delta_2(P)$, which correspond to the difference between the extremum and the second extremum on lines \ref{line:maxindex} and \ref{line:remindex}, respectively.

    With all of the notation in place, we define
    \begin{align*}
        \epsilon_0 &= \min_{P\in S_1 \setminus \{P_0\}} \frac{\max(\bs{d}(P))}{\|\bs{V}(\bs{I} - {\bs{V}_{P*}}^T({\bs{V}_{P*}}^{T})^\dagger)\bs{V}^T\|_2},\\
        \epsilon_1 &= \min_{P\in S_1\setminus \{P_0\}} \frac{\delta_1(P)/2}{\|\bs{V}(\bs{I} - {\bs{V}_{P*}}^T({\bs{V}_{P*}}^{T})^\dagger)\bs{V}^T\|_2},\\
        \epsilon_2 &= \min_{P\in S_2} \frac{\min_{i\in P}|s(P)_{i}|}{\|({\bs{V}_{P*}}^{T})^\dagger\bs{V}^T\|_2},\\
        \epsilon_3 &= \min_{P \in S_2}\,\min_{i\in Q(P)} \frac{w(P)_i-s(P)_i}{2\left(\|(\bs{V}_{P'*}^T)^\dagger\bs{V}^T\|_\infty+\|(\bs{V}_{P*}^T)^\dagger\bs{V}^T\|_\infty\right)},\\
        \epsilon_4 &= \min_{P\in S_2}\,\min_{i\in Q(P)} \frac{\delta_2(P)\left(w(P)_i-s(P)_i\right)^2}{4\left(-s(P)_i\|(\bs{V}_{P'*}^T)^\dagger\bs{V}^T\|_\infty + w(P)_i\|(\bs{V}_{P*}^T)^\dagger\bs{V}^T\|_\infty\right)},\\
        \epsilon &= \min\{\epsilon_0,\epsilon_1,\epsilon_2,\epsilon_3,\epsilon_4\}.
    \end{align*}
    We set $\delta_2(P) = +\infty$ when $|Q(P)| \leq 1$, and take a minimum over an empty set to be $+\infty$. By \Cref{assum:nnls}, $\delta_1(P) > 0$ and $\delta_2(P) > 0$ for all $P$, and since $S_1$ and $S_2$ are finite, $\epsilon > 0$.
    Now assume a perturbed measure with weights $\widetilde{\bs{w}} = \bs{w} + \Delta\bs{w}\in\R^M$ satisfying $\|\Delta\bs{w}\|_2<\epsilon$. For any $P\in S_1\setminus \{P_0\}$,
    \begin{align*}
        \widetilde{\bs{d}}(P) &= \bs{V}(\widetilde{\bs{\eta}} - {\bs{V}_{P*}}^T\widetilde{\bs{s}}(P)_P) = \bs{V}(\bs{V}^T\widetilde{\bs{w}} - {\bs{V}_{P*}}^T\widetilde{\bs{s}}(P)_P)\\
        &= \bs{V}(\bs{V}^T(\bs{w} + \Delta\bs{w}) - {\bs{V}_{P*}}^T(\bs{s}(P)_P + ({\bs{V}_{P*}}^{T})^\dagger \bs{V}^T\Delta{\bs{w}}))\\
        &= \bs{V}(\bs{\eta} - {\bs{V}_{P*}}^T\bs{s}(P)_P + \bs{V}^T(\Delta\bs{w}) - {\bs{V}_{P*}}^T(({\bs{V}_{P*}}^{T})^\dagger \bs{V}^T\Delta{\bs{w}}))\\
        &= \bs{d}(P) + \bs{V}(\bs{I} - {\bs{V}_{P*}}^T({\bs{V}_{P*}}^{T})^\dagger)\bs{V}^T\Delta{\bs{w}},
    \end{align*}
    which implies,
    \begin{align*}
        \|\widetilde{\bs{d}}(P) - \bs{d}(P)\|_2 &= \|\bs{V}(\bs{I} - {\bs{V}_{P*}}^T({\bs{V}_{P*}}^{T})^\dagger)\bs{V}^T\Delta{\bs{w}}\|_2\\
        &< \|\bs{V}(\bs{I} - {\bs{V}_{P*}}^T({\bs{V}_{P*}}^{T})^\dagger)\bs{V}^T\|_2\epsilon \leq \max(\bs{d}(P)),
    \end{align*}
    and this last strict inequality ensures that $\max(\widetilde{\bs{d}}) >0$, so that for any iteration of the unperturbed algorithm when the outer loop is triggered on line \ref{line:nnlsouterloop}, the perturbed algorithm also has this condition triggered. Building on the inequality above, we have,
    \begin{align*}
      \|\widetilde{\bs{d}}(P) - \bs{d}(P)\|_\infty \leq \|\widetilde{\bs{d}}(P) - \bs{d}(P)\|_2 < \|\bs{V}(\bs{I} - {\bs{V}_{P*}}^T({\bs{V}_{P*}}^{T})^\dagger)\bs{V}^T\|_2\epsilon \leq \frac{\delta_1(P)}{2},
    \end{align*}
    for all $P\in S_1\setminus \{P_0\}$ ensuring that the maximization problem on line \ref{line:maxindex} has the same solution in both the perturbed and unperturbed algorithms. To ensure $Q(P)$ defined on line \ref{line:Qdef} is the same in both perturbed and unperturbed algorithms, we compute,
    \begin{align*}
      \|\widetilde{\bs{s}}(P)_P - \bs{s}(P)_P\|_2 &= \| ({\bs{V}_{P*}}^{T})^\dagger\widetilde{\bs{\eta}} - \bs{s}(P)_P \|_2 \\
      &= \| \bs{s}(P)_P + ({\bs{V}_{P*}}^{T})^\dagger \bs{V}^T\Delta{\bs{w}} - \bs{s}(P)_P \|_2 \\
      &< \|({\bs{V}_{P*}}^{T})^\dagger \bs{V}^T\|_2\epsilon \stackrel{\epsilon \leq \epsilon_2}{\leq} \min_{i\in P}|s(P)_{i}|,
    \end{align*}
    for all $P\in S_2$ ensuring that every entry of $\widetilde{\bs{s}}(P)_P$ has the same sign as the corresponding entry of $\bs{s}(P)_P$, so $\widetilde{Q}(P) = Q(P)$.
    Finally, we ensure that the minimization on line \ref{line:remindex} identifies the same index in both perturbed and unperturbed cases whenever $|Q(P)|\geq2$. 
    Both algorithms zero the same entries, so $\widetilde{w}(P)_i - w(P)_i = a_i \coloneqq ((\bs{V}_{P'*}^T)^\dagger\bs{V}^T\Delta\bs{w})_i$ for $i \in P\cap P'$, and $\widetilde{w}(P)_m = w(P)_m = 0$. Similarly, $\widetilde{s}(P)_i - s(P)_i = b_i \coloneqq ((\bs{V}_{P*}^T)^\dagger\bs{V}^T\Delta\bs{w})_i$. If the index $m$ added in the current outer iteration satisfies $m \in Q(P)$, its ratio on line \ref{line:remindex} is $0$, and it is removed. For $i \in Q(P)\backslash\{m\}$, we have $s(P)_i < 0 < w(P)_i$, $|a_i| \leq \|(\bs{V}_{P'*}^T)^\dagger\bs{V}^T\|_\infty\|\Delta\bs{w}\|_2$ and $|b_i| \leq \|(\bs{V}_{P*}^T)^\dagger\bs{V}^T\|_\infty\|\Delta\bs{w}\|_2$, so $\epsilon \leq \epsilon_3$ gives $|b_i - a_i| < (w(P)_i - s(P)_i)/2$. Then,
    \begin{align*}
        &\left|\frac{-\widetilde{w}(P)_i}{\widetilde{s}(P)_i - \widetilde{w}(P)_i} - \frac{-{w}(P)_i}{{s}(P)_i - {w}(P)_i}\right|\\
        &= \left| \frac{-a_i\, s(P)_i + w(P)_i\, b_i}{\big(s(P)_i-w(P)_i+b_i-a_i\big)\big(s(P)_i-w(P)_i\big)} \right|\\
        &\leq \frac{2\left|-a_i\, s(P)_i + w(P)_i\, b_i\right|}{(s(P)_i-w(P)_i)^2} \\
        &< \frac{2\left(-s(P)_i\|(\bs{V}_{P'*}^T)^\dagger\bs{V}^T\|_\infty+w(P)_i\|(\bs{V}_{P*}^T)^\dagger\bs{V}^T\|_\infty\right)}{(s(P)_i-w(P)_i)^2}\,\epsilon \\
        &\stackrel{\epsilon\leq\epsilon_4}{\leq} \frac{\delta_2(P)}{2},
    \end{align*}
    which is half the optimality gap, so that the minimum index obtained in line \ref{line:remindex} is unchanged.
    
    Finally, due to exactness of the least squares solution with $|P_0|=N$, we have that $\widetilde{\bs{d}}(P_0) = \bs{0}$, exiting the outer loop at the same time as the unperturbed problem, and the solution satisfies, with $\bs u = \bs s(P_0)_{P_0}$ and $\widetilde{\bs u} = \widetilde{\bs s}(P_0)_{P_0}$,
    \[\|\widetilde{\bs{u}} - \bs{u}\|_2 =  \|({\bs{V}_{P_0*}}^{T})^\dagger \bs{V}^T\Delta{\bs{w}}\|_2 < \|({\bs{V}_{P_0*}}^{T})^\dagger \bs{V}^T\|_2\epsilon.\]
\end{proof}
The same argument applies to the full Lawson-Hanson method \cite{lawson_1995}, whose inner loop additionally moves the iterate, $\bs{w} \gets \bs{w} + \alpha(\bs{s} - \bs{w})$, with $\alpha$ the minimum ratio on line \ref{line:remindex}. The only change is that $\bs{w}(P)$ is no longer linear in $\Delta\bs{w}$: along a fixed sequence of index sets, it is a composition of finitely many rational functions of $\bs{\eta}$ whose denominators $w(P)_i - s(P)_i$ are bounded away from zero. The same computation as above bounds the perturbation of each step size, $\alpha$. Inductively, this yields constants $L(P)$, depending only on the unperturbed run, such that $\|\widetilde{\bs{w}}(P) - \bs{w}(P)\|_\infty \leq L(P)\|\Delta\bs{w}\|_2$. Replacing $\|(\bs{V}_{P'*}^T)^\dagger\bs{V}^T\|_\infty$ by $L(P)$ in $\epsilon_3$ and $\epsilon_4$ then gives the result.

\end{document}